\documentclass[11pt]{article}
\usepackage{amsmath, amssymb}
\usepackage[french]{babel}
\usepackage[dvips]{graphics}
\usepackage{graphicx}
\usepackage{psfrag}
\usepackage{graphics}
\usepackage{verbatim}

\usepackage[french]{babel}

\title{Combinatorial aspects of matrix models}

\author{Alice Guionnet\thanks{Ecole Normale
 Sup\'erieure de Lyon,
Unit\'e de Math\'ematiques pures et appliqu\'ees,
UMR 5669,
46 All\'ee d'Italie, 
69364 Lyon Cedex 07, France. E-mail: aguionne@umpa.ens-lyon.fr.},
Edouard Maurel-Segala\thanks{Ecole Normale Sup\'erieure de Lyon,
Unit\'e de Math\'ematiques pures et appliqu\'ees,
UMR 5669,
46 All\'ee d'Italie, 
69364 Lyon Cedex 07, France. E-mail: emaurel@umpa.ens-lyon.fr}}

\newtheorem{prop}{Proposition}[section]
\newtheorem{theo}[prop]{Theorem}

\newtheorem{lem}[prop]{Lemma}

\newtheorem{pr}[prop]{Property}
\newtheorem{cor}[prop]{Corollary}

\newenvironment{dem}{\textbf{Proof.}\par}
{\begin{flushright}$\Box$\end{flushright}}

\def\E{\mathbb E}
\def\C{\mathbb C}
\def\e{\epsilon}
\newcommand{\RR}{\mathbb{R}}
\newcommand{\CC}{\mathbb{C}}
\newcommand{\NN}{\mathbb{N}}

\def\N{{\mathbb N}}

\def\tr{{\mbox{tr}}}
\def\trn{{\tr_N}}

\def\nn{\noindent}
\def\oo{\overline}

\def\bA{{\bf A}}

\def\bA{{{\bf A}}}

\def\a{\alpha}
\def\b{\beta}
\def\d{\delta}
\def\e{\epsilon}
\def\g{\gamma}

\def\l{\lambda}

\def\s{\sigma}

\def\D{\Delta}
\def\G{\Gamma}

\def\ra{\rightarrow}

\def\Aa{{\cal A}}

\def\Ca{{\cal C}}
\def\Da{{\cal D}}

\def\Ha{{\cal H}}

\def\La{{\cal L}}
\def\Ma{{\cal M}}
\def\MM{{\cal M}}
\def\DD{{\cal D}}
\def\JJ{{\cal J}}
\def\II{{\cal I}}

\def\mun{{\hat\mu^N}}

\def\lbc{\lbrace}

\def\nn{\noindent}
\def\part{\partial}

\def\ot{\otimes}

\def\ts{\times}
\def\Pa{{\mathcal{P}}}

\def\R{{\mathbb R}}

\def\cxm{\C\langle X_1,\cdots, X_m\rangle}

\def\bX{{\bf X}}

\def\bM{{\bf M}}

\begin{document}

\renewcommand{\refname}{\Large{References}}

\maketitle
\centerline{\bf Abstract}
We show that under reasonably general
assumptions, the first order asymptotics 
of the free energy of matrix models 
are generating functions for colored planar
maps. This is based on the fact that
solutions of the Schwinger-Dyson equations
are, by nature, generating functions
for enumerating planar maps, a remark which bypasses
the use of 
Gaussian calculus.

\medskip

\nn
{\it Keywords :} Random matrices, non-commutative measure, map enumeration.

\medskip

\nn
{\it Mathematics Subject of Classification :}   15A52, 46L50, 05C30.

\nn
\section{Introduction}
It has long been used in combinatorics and  physics 
 that moments of Gaussian matrices 
have a valuable combinatorial interpretation. 
The first result
in this direction was due to Wigner \cite{Wig}
who proved that 
the trace of even moments of 
a $N\ts N$ Hermitian matrix $A$ with 
i.i.d centered entries with covariance 
$N^{-1}$ converge as $N$ goes to infinity
towards the Catalan numbers which enumerate 
  non crossing partitions. 
If one restricts to Gaussian  entries,
that is matrices following the law $\mu_N$
of the {\bf GUE}  which is the law on the set $ \Ha_N$
of $N\ts N$ Hermitian matrices 
with density
$$\mu_N(dA)=\frac{1}{Z_N}1_{A\in \Ha_N} e^{-\frac N 2 \tr(A^2)}
\prod_{1\le i\le j \le N}d\Re e(A_{ij}) 
\prod_{1\le i< j \le N}d\Im m(A_{ij}),$$
it occurs that the corrections 
to this convergence count
graphs which can be embedded 
on surface of higher genus, a fact which was used by  Harer and Zagier \cite{harerzagier}.
This enumerative property was fully developed
after 't Hooft, who  considered generating functions
of such moments.

For instance,
c.f  Zvonkin \cite{zvon}, we have the { formal } expansion
 $$F_N(tx^4)=\frac{1}{N^2}\log \int e^{-N t\tr(A^4)}d\mu_N(A)
=\sum_{k\ge 1}\sum_{g\ge 0}
\frac{(-t)^k}{ k! N^{2g}} C(k,g)$$ 
with
$$\begin{array}{rl}
C(k,g)=\mbox{Card}\{& \hspace{-0.4cm}\textrm{ maps with genus } g\\
& \hspace{-0.4cm}\textrm{ with $k$ stars of valence $4$}\}
\end{array}$$

Here, maps are connected oriented diagrams which
can be embedded into a surface of genus 
$g$ in such a way that edges do not cross
and the faces of the graph (which are defined by 
following the boundary of the graph) are homeomorphic
to a disc. The counting is
done up to equivalent classes, i.e. up to homeomorphism.
Let us stress that the above equality is  
only formal and 
should be understood in the sense
that all the derivatives at the origin on both
sides of the equality match, it means that, for all $k\in\N$,
$$ (-1)^k\partial_t^k F_N(tx^4)|_{t=0}
=\sum_{g\ge 0}
\frac{1}{N^{2g}} C(k,g)$$
which can be proved thanks to  Wick's formula (note above that
the sum is in fact finite).

Such expansions can be generalized to arbitrary 
polynomial functions (to enumerate maps 
with vertices of different degrees)
and to 
several-matrices integrals which allow
to enumerate colored maps; if  $V$ 
 is a self-adjoint    polynomial
of $m$ non-commutative variables,
$$F_N(V)=\frac{1}{N^2}\log \int e^{-N\tr(V(A_1,\cdots,A_m))} d\mu_N(A_1)
\cdots d\mu_N(A_m)$$
expands, when $V=\sum_{i=1}^n t_i (q_i+q_i^*)$ with some monomials $q_i$
and real parameters $t_i$, with $q_i^*$ being
the adjoint of $q_i$ (see section \ref{sd})
into an enumeration of colored maps.

Our aim is to look beyond this formal work and gives
a rigorous proof of this expansion.

In the case of one matrix integrals, this problem
is quite well understood
at any level of the expansion
and for  any reasonable 
potentials $V$ (see \cite{APS} and \cite{EML}
for instance).

Several matrix models are much harder. In the physics 
literature, the focus 
is mostly on a few specific 
integrals;
we refer the interested reader to the reviews 
 \cite{FGZ,GPW}. In the mathematical literature, 
fewer matrix integrals could be analyzed and only their first
order asymptotics could be derived (see 
Mehta et al. \cite{M-M,Me2} and Guionnet \cite{GCMP,GM}).
In free probability, even the problem of the existence of the
free energy is wide open for reasonably general potential $V$
and its solution would have important consequences.
In combinatorics, another road was opened by 
Bousquet-Melou and Schaeffer \cite{BM-S},
following the ideas of Tutte \cite{Tu},
by using directly bijection between 
maps and well labeled trees.
 
In this paper, we shall rather focus on an even 
more interesting quantity 
than the free energy, namely,  the limiting empirical 
distribution of matrices; for $A_1,\cdots,
A_m\in\Ha_N^m$, it is defined
as the linear form on  the set $\cxm$ of polynomials
of $m$ non-commutative variables so that
$$\mun_{A_1,\cdots, A_m}(P)=\frac{1}{N} \tr(P(A_1,\cdots, A_m)).$$
 Let
$\mu^N_V$ be the law on $\Ha_N^m$ given by
$$\mu^N_V(dA_1,\cdots, dA_m)=e^{-N^2 F_N(V)}
e^{-
N\tr(V(A_1,\cdots,A_m))}
\prod_{i=1}^m d\mu_N(A_i)$$
with $F_N(V)$ as above. 

For the following, we take a potential $V=V_{\oo{t}}=\sum_{i=1}^n t_i (q_i+q_i^*)$.
Then we proceed in two steps to relate the first asymptotic of $\mu^N_V$
to the enumeration of planar graphs.
\goodbreak
$\bullet$
First, we study the solution $\tau\in\cxm^*$ of the so-called Schwinger-Dyson equations {\bf SD[V]}:

$$\tau_{V}((X_i+\Da_i V)P)=\tau_{V}\otimes\tau_{V}(D_i P)$$

for all $P\in \cxm$ and $i\in\{1,\cdots,m\}$. 
Here, $D_i$ and $\Da_i$ are respectively
the non-commutative derivative and the cyclic derivative
with respect to the $i^{th}$ variable (see paragraph \ref{deri}).
More precisely we're interested in properties of unicity and existence
for the solution to this equation.

Moreover, let us associate to $(X_i)_{1\le i\le m}$ 
$m$ branches of different colors, 
and to a monomial $q(\bX)=X_{i_1}\cdots X_{i_p}$
a star with $p$ colored 
branches by ordering clockwise the branches corresponding
to $X_{i_1},\cdots, X_{i_p}$.
Such a star is said to be of type $q$. Note that it has 
a distinguished branch, the first one,  $X_{i_1}$,
and its branches are oriented by the above
 clockwise order (one should imagine the star to be 
fat, each branch made of two parallel segments
 which have opposite orientation, the whole orientation being given
by the clockwise order).
This defines a bijection between non-commutative
monomial and oriented stars with colored branches and one
distinguished branch.
(see a precise description of the planar maps we enumerate in
subsection \ref{combinatorics}). 

We can now relate Schwinger Dyson's equation and maps
enumeration:

\begin{theo}\label{theo1}
Let $R>2$,
then there exists an open  neighborhood $U\subset{\mathbb R}^n$ 
of the origin
(a ball of positive radius)
  such that:
\begin{itemize}
\item
For ${\oo{t}}\in U$, there exists a unique $\tau_{\oo{t}}\in\cxm^*$ 
which is a solution to {\bf SD[$V_{\oo{t}}$]} and such that for all $p$,
for all $i_1,\cdots,i_p$ in $[|1,m|]$,
$\tau_{\oo{t}}(X_{i_1}\cdots X_{i_1})\leq R^p$
\item
For all $P$ monomial in $\cxm$,
${\oo{t}}\ra \tau_{\oo{t}}(P)$ is analytic on $U$ and for all $k_1,\cdots,k_n$ integers,
$(-1)^{\Sigma k_i}\partial_{t_1}^{k_1}\cdots\partial_{t_n}^{k_n}
\tau_{\oo{t}}(P)|_{\oo{t}=0}$ is the number of maps with $k_i$ stars of type $q_i$ or $q_i^*$ and one of type $P$.
\end{itemize}
\end{theo}

$\bullet$
Then, we shall see (see section \ref{mm})
that, under some appropriate assumptions on $V$,
$\mun_{A_1,\cdots, A_m}$ converges under
$\mu^N_V$ towards a solution $\tau_V$ to the  
Schwinger-Dyson equations {\bf SD[$V$]}.

First if $V$ is sufficiently large, 
the limit will solve a weak form of the Schwinger-Dyson equation (see section \ref{limitpoint}).
Then we will consider convex potential $V$ (see section \ref{convexsec}) 
for which we have:

\begin{theo}
Let $U_a$ be the set of $t_i$'s
for which $V_{\oo{t}}$ is convex, then there exists $\e>0$ 
such that for  $(t_i)_{1\le i\le n} \in U_a\cap B(0,\e)$, $\mu^N_{V_{\oo{t}}}[\mun_\bA]$ converges to
the unique solution to {\bf SD[$V_{\oo{t}}$]} as described in theorem \ref{theo1}
\end{theo}

Hence, $(\tau_{V_{\oo{t}}})_{|t|\le \e}$
 are {generating states } 
for the enumeration 
of colored planar maps and  Schwinger Dyson's equations
can be viewed as the generating differential 
equations to enumerate  
colored planar maps. This is due
to the fact that the action of the derivatives
$D_i$ and $\Da_i$ on monomials, under the above bijection
between stars and monomials, produces natural operations
on planar maps.  The fact that they are related 
with matrix integrals can be used 
to actually study these equations
and obtain informations
about their solutions, henceforth
solving the related combinatorial
problem. It further should  give information over higher
genus maps, a problem that we shall tackle in a forthcoming paper.

Coming back to the
free energy of matrix models, 
we conclude (see Theorem \ref{convfe})
that when the empirical distribution
of matrices converges towards the
solution to Schwinger-Dyson's 
equations, 
the free energy  is also a generating function
of the associated planar maps. In section \ref{div}, we give another
sufficient condition for this convergence
to hold. More precisely, we argue that the whole machinery works
for a general polynomial $V_{\oo{t}}$ with $\oo{t}$ sufficiently small
if we add a well chosen cut-off.

Finally, we will apply these results to
the study of Voiculescu's microstates entropy
and summarize  direct applications to the enumeration of 
planar maps associated with some matrix models.

The results of this paper
are clearly 
known, at least at a subconscious
level, by physicists but we could 
not find any proper reference on the subject. 
On a mathematical level, it is rather elementary
and we  hope it will demystify this
interesting field of physics to  mathematicians, 
or at least to probabilists.

\section{Schwinger Dyson's equations and combinatorics}\label{sd}
\subsection{Tracial states}
Let $\cxm$ be the set of polynomial functions in $m$ self-adjoint
non-commutative variables.
We endow $\cxm$ with the involution given for all $z\in\C$, all $i_1,\cdots,i_p\in\{1,\cdots,m\}$
and all $p\in\N$, by
$$(zX_{i_1}\cdots X_{i_p})^*=\bar z X_{i_p}\cdots X_{i_1}.$$
We will say that $P$ in $\cxm$ is self-adjoint if $P^*=P$.

For any $R>0$, completing 
$\cxm$ for  the
norm
$$||P||_R=\sup_{\Aa\ C^*-\mbox{algebra}}
\sup_{\genfrac{}{}{0pt}{}{a_1,..,a_m\in \Aa,}{a_i=a_i^*\|a_i\|_\Aa\le R}}
 \|P(a_1,\cdots,a_m)\|_\Aa$$
 produces a  $\C^*$-algebra $\cxm_R=(\cxm, ||.||_R,*)$.

We let $\cxm^*$ be the set of
real linear forms on $\cxm$ (i.e linear forms such that
$\tau(a^*)=\overline{\tau(a)}$), and denote 
$\cxm^*_R$ the subset of $\cxm^*$ 
of continuous forms with respect to
the norm $||.||_R$, i.e the topological dual of $\cxm_R$.

We let $\Ma^m$ be the set 
of laws of $m$ bounded self-adjoint non-commutative
variables, that is the subset of  elements
$\tau$ of $\cxm^*$ such that

\begin{equation}\label{tracial}
\tau(PP^*)\ge 0,\quad \tau(PQ)=\tau(QP)\quad\forall P,Q\in
 \C\langle X_1,\cdots,X_m\rangle, \quad 
\tau(I)=1.
\end{equation}
It is not hard to see that for any $R<\infty$, $\Ma^m_R=\cxm_R^*\cap \Ma^m$
is a compact metric space for the weak*-topology
by Banach-Alaoglu theorem. Elements of $\Ma^m=\cup_{R\ge 0}\Ma^m_R$
are said to be compactly supported, by analogy with the 
case $m=1$ where they are indeed
compactly supported probability measures. A family $(\tau_t)_{t\in I}$
of elements of $\Ma^m_R$ for some $R<\infty$ is
said to be uniformly compactly supported.

To deal with  variables which do not
have all their moments, we eventually can change 
the set of test functions and, following \cite{CDG1},
 consider
instead of $\cxm$ the  complex vector space $\Ca^m_{st}(\C)$ generated by
the Stieljes functionals
\begin{equation}\label{defst}
ST^m(\C)=\{ 
\prod_{1\le i\le p}^{\ra}(z_i-\sum_{k=1}^m\a_i^k
\bX_k)^{-1}
;\quad  z_i\in\C\backslash\R, \a_i^k\in \R, p\in\N\}
\end{equation}
where $\prod^\ra$ is the non-commutative product.
We can give to $ST^m(\C)$ an involution 
and a norm 
$$||F||_\infty=\sup_{\Aa C^*-\mbox{algebra}}\sup_{a_i=a_i^*\in\Aa}||F(a_1,\cdots,a_m)||_\infty$$
where the supremum is taken eventually
on unbounded operators affiliated with $\Aa$,
which turns it into a $C^*$-algebra.
We will denote  $\Ca_{st}^m(\R)=\{G=F+F^*, F\in \Ca_{st}^m(\C)\}
 $.
We will let $\Ma^m_{ST}$ be the set of linear forms on $\Ca_{st}^m(\C)$
which satisfy (\ref{tracial}) (but with functions 
of $\Ca_{st}^m(\C)$ instead of $\cxm$). If one equips $\Ma^m_{ST}$ 
with its weak topology, then $\Ma^m_{ST}$ 
is a compact metric space (see \cite{CDG1}).

\subsection{Non-commutative derivatives}\label{deri}
We let $D_i:\cxm\ra\cxm\otimes\cxm$ be the non-commutative 
derivative with respect to $X_i$ given
by the  Leibnitz rule 
$$D_i(PQ)=D_iP\ts (1\otimes Q)+(P\otimes 1 )\ts D_iQ$$
for any $P,Q\in\cxm$
and the condition
$$D_iX_j=1_{i=j} 1\otimes 1.$$
In other words, if $P$ is a non-commutative monomial
$$D_iP=\sum_{P=P_1X_iP_2} P_1\otimes P_2$$
where the sum runs over 
over all possible decomposition of $P$ as $P_1X_iP_2$.
This definition can be extended to
$\Ca^m_{st}(\C)$ by keeping the above Leibnitz rule
(but with $P,Q$ in $\Ca^m_{st}(\C)$) and 
$$D_i (z_i-\sum_{k=1}^m\a_k
\bX_k)^{-1}=\a_i(z_i-\sum_{k=1}^m \a_k \bX_k)^{-1}\otimes
(z_i-\sum_{k=1}^m\a_k
\bX_k)^{-1}.$$
We also define the cyclic derivative $\Da_i$
as follows. Let $m:\cxm\otimes\cxm\ra\cxm$
(resp. $\Ca_{st}^m(\C)\otimes\Ca_{st}^m(\C)\ra\Ca_{st}^m(\C)$)
be defined by $m(P\otimes Q)=QP.$
Then, we set $$\Da_i=m\circ D_i.$$
If $P$ is a non-commutative monomial,
we have
$$\Da_i P=\sum_{P=P_1 X_iP_2} P_2 P_1.$$
\subsection{Schwinger-Dyson's equation}
Let $V$ be self-adjoint and 
 consider the following equation on 
$\cxm^*$; we say that $\tau\in\cxm^*$
satisfies the Schwinger-Dyson equation
with potential $V$, denoted in short {\bf SD[V]},
if and only if for all $i\in\{1,\cdots,m\}$ and $P\in\cxm$,

$$\tau(I)=1,\quad\tau\otimes \tau(D_i P)=
\tau( (\Da_i V+X_i)P)\quad\quad\quad {\bf SD[V]}$$

These equations are called Schwinger-Dyson's equations in physics,
but
in free probability, one would rather say that
the conjugate variable (or alternatively the 
non-commutative Hilbert transform) $D^*_i1$ under $\tau$
is equal to $X_i+\Da_i V$ for all $i\in\{1,\cdots,m\}$.

\subsection{Uniqueness of the solutions
to Schwinger-Dyson's equations for small parameters}
In this paper, we shall restrict ourselves 
to non-oscillatory integrals, that is to
the case where $\tr(V(X_1,\cdots,X_m))$
is real for any $m$-uple 
of Hermitian matrices. In other words, $$\tr(V(X_1,\cdots,X_m))
=\tr(V^*(X_1,\cdots,X_m))=\tr( 2^{-1} (V+V^*)(X_1,\cdots,X_m))$$
for any $m$-uple 
of Hermitian matrices. Thus, we shall 
 assume that 
$$V(X_1,\cdots,X_m)=V_{\oo{t}}(X_1,\cdots,X_m)=
\sum_{i=1}^n t_i (q_i(X_1,\cdots,X_m)+q_i^*(X_1,\cdots,X_m))$$
where the $q_i$'s are monomial functions of $m$ non-commutative indeterminates  and  $\oo{t}=
(t_1,\cdots, t_n)$ 
are real parameters.

In this paragraph, we shall 
consider solutions to {\bf SD[$V_{\oo{t}}$]}
which satisfy a compactness condition
that we shall discuss in the following subsections.
Let  $R\in {\mathbb R}^+$
(We will always assume $R\ge 1$ without loss
of generality).

{\bf (H(R))}{\it  An element $\tau\in \cxm^*$
satisfies {\bf (H(R))} iff   for all $k\in\N$,
$$\max_{1\le i_1,\cdots,i_k\le m}
| \tau(X_{i_1}\cdots X_{i_k})| \le R^k.$$ }

In the sequel, we denote $D$ the degree of $V$, that is 
the maximal degree of 
the $q_i's$; $q_i(X)=X_{j_1^i}\cdots X_{j_{d_i}^i}$
with, for $1\le i\le n$,  $\mbox{deg}(q_i)=:
d_i\le D$ and equality holds 
for some $i$.

The main result of this paragraph  is 

\begin{theo}\label{unique}
For all $\oo{t}\in\R^n$, there exists $A(\oo{t})=A(|t|)\in\R^+$
with $|t|=\max_{1\le i\le n} |t_i|$, $A(|t|)$ goes to infinity
when $|t|$ goes to zero, so that for $R\le A(|t|)$,
there exists at most one solution 
$\tau_{\oo{t}}$ to {\bf SD[$V_{\oo{t}}$]} which satisfies 
{\bf (H(R))}.
\end{theo}

\nn
{\bf Remark: }
Note here that  it could be believed at first sight
that the solutions to {\bf SD[V]}
are not unique since they depend 
on the trace of high moments $\tau(q_jP)$.
However, our compactness assumption {\bf (H(R))}
gives uniqueness because it forces the
solution to be in a small neighborhood 
of the law $\tau_0=\s^m$ of $m$ free semi-circular variables,
so that  perturbation analysis
applies. We shall see in Theorem \ref{existanaltheo}
that this solution is actually the one
which is related with the enumeration of maps. 

\nn
\begin{dem}
Let us assume we have two solutions $ \tau$ and $ \tau'$.
Then, by the equation {\bf SD[V]},
for any monomial function $P$ of degree $l-1$, for
$i\in\{1,\cdots,m\}$,
$$( \tau- \tau')(X_iP)=(( \tau-\tau')\otimes\tau)(D_i P)+(\tau'\otimes(\tau-\tau'))(D_i P) -( \tau- \tau')(\Da_i V P)$$
Hence, if we let for $l\in\N$
$$\D_l( \tau, \tau')=
\sup_{\mbox{monomial } P\mbox{  of degree   l}}
|
 \tau(P)- \tau'(P)|$$
we get, since if $P$ is of degree $l-1$, $$D_iP=\sum_{k=0}^{l-2}
p^1_k\otimes p^2_{l-2-k}$$
where  $p_k^i$, $i=1,2$ are
monomial of degree $k$ or the null monomial, 
and $\Da_i V$ is a finite sum of
monomials of degree smaller than $D-1$,
$$\D_l( \tau, \tau')=
\max_{P \mbox{  of degree   } l-1}\max_{1\le i\le m}\{ | \tau(X_iP)- \tau'(X_iP)|\}$$
$$
 \le 2\sum_{k=0}^{l-2} \D_k( \tau, \tau')R^{l-2-k}
+ C|t| \sum_{p=0}^{D-1}  \D_{l+p-1}( \tau, \tau')$$
with a finite constant $C$ (which depends on $n$ only).
For $\gamma>0$, we set
 $$d_\gamma({ \tau, \tau'})=\sum_{l\ge 0} \gamma^l
\D_l( \tau, \tau').$$
Note that under {\bf (H(R))},
this sum is finite for $\g<(R)^{-1}$.
Summing the two sides  of the above inequality 
times $\gamma^l$
we arrive at
$$d_\gamma({ \tau, \tau'})
\le 2\gamma^2  (1-\gamma R )^{-1}d_\gamma({ \tau, \tau'})
+C|t| \sum_{p=0}^{D-1}\gamma^{-p+1} d_\gamma({ \tau, \tau'}).
$$
We finally conclude that if $(R,|t|)$ are
 small enough so that we can choose 
$\gamma\in (0, R^{-1})$ so that
$$2\gamma^2 (1-\gamma R )^{-1}
+C|t|\sum_{p=0}^{D-1}\gamma^{-p+1} <1$$
$d_\gamma({ \tau, \tau'})=0$ and so $ \tau= \tau'$
and 
we have at most one solution. Taking $\gamma=(2R)^{-1}$
shows that this is possible provided 
$$\frac{1}{4R^2} +C|t|\sum_{p=0}^{D-1} (2R)^{p-1}<1$$
so that when $|t|$ goes to zero, we see that we 
need $R$ to be smaller than $A(|t|)$ of order
$|t|^{-\frac{1}{D-2}}$.

\end{dem}

\subsection{Combinatorics}\label{combinatorics}
In this paragraph we describe the combinatorial objects 
we are considering.

To describe the enumeration we have to deal
with, let us associate a colored star 
to any monomial. We associate to each
 $i\in \{1,\cdots,m\}$ a different
color. Then, we define a
bijection between oriented  branch-colored 
stars with a distinguished branch 
and non-commutative monomials as follows.
For any $i\in \{1,\cdots,m\}$,
we associate to  $X_i$  a branch of color $i$.
 We shall say that
a star is of type $q(X_1,\cdots,X_m)=X_{i_1}\cdots X_{i_l}$
if it is a star with $l$ branches which we color 
clockwise ; the first branch will be of color $i_1$, 
the second of color $i_2$
... etc ... until the $l^{th}$ branch is colored with color $i_l$. 
Note that this star possesses a
distinguished branch, the one corresponding
to $X_{i_1}$, and an orientation,
corresponding to the clockwise
order.
By convention,
the star of type $q=1$ is simply 
a point.

\begin{figure}[ht!]
\psfrag{x1}{\begin{LARGE}{$X_{1}$}\end{LARGE}}
\psfrag{x2}{\begin{LARGE}{$X_{2}$}\end{LARGE}}
\psfrag{m1}{\begin{Large}{marks the first branch}\end{Large}}
\psfrag{m2}{\begin{Large}{indicates the orientation}\end{Large}}
\begin{center} \resizebox{!}{5cm}{\includegraphics{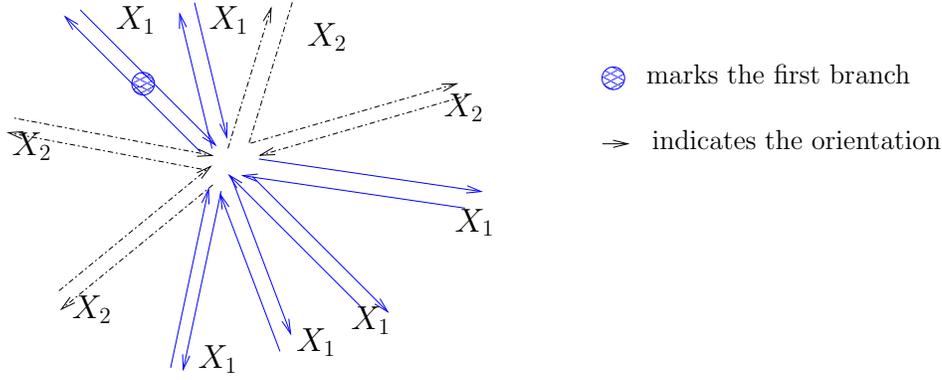}}
\end{center}
\caption{The star of type $q(X)=X_1^2X_2^2X_1^4X_2^2$}
\end{figure}

\par
A planar map is a connected graph embedded into the sphere
with colored stars, each branch is glued with exactly one branch 
of the same color and the
edges obtained in this way do not cross each other. Hence, branches
are thought as half edges.  Maps are only considered up to an homeomorphism of
the sphere. Now we will be interested in enumerating maps 
with a fixed set of stars, we define for $q_i$ a 
family of non necessarily 
distinct monomials and $k_i$ a family of integers:
$$\MM_0((q_1,k_1),(q_2,k_2),\ldots,(q_n,k_n))=
\sharp\{\textrm{planar maps build with $k_i$ stars of type $q_i$}\}.$$
We denote in short $\MM_0(P,(q_1,k_1),\ldots,(q_n,k_n))=
\MM_0((P,1),(q_1,k_1),\ldots,(q_n,k_n))$.
In that set, each star is labeled and has a marked branch
 (which correspond to its first variable) so that for
example $\MM_0((X^4,2))=36$.
Now what we are really interested in is enumerating maps with a fixed number of stars of type $q$ or $q^*$ so that we define:

$$\MM((q_1,k_1),\ldots,(q_n,k_n))=\sum_{\genfrac{}{}{0pt}{}{1\leq p_i\leq k_i,}{1\leq i\leq n}}\prod_{i=1}^nC^{p_i}_{k_i}\MM_0((q_1,p_1),(q_1^*,k_1-p_1),\ldots,(q_n,k_n),(q_n^*,k_n-p_n))$$
and $\MM(P,(q_1,k_1),\ldots,(q_n,k_n))=\MM((P,1),(q_1,k_1),\ldots,(q_n,k_n)).$

This quantity enumerates the number of ways
 to build a map on stars of fixed types up to the symmetry induced on stars by the operand $*$.

\par
Due to the fact that everything is labeled,
 we enumerate lots of very similar objects.
A way to avoid this problem is to look at the maps as 
they are enumerated by
combinatoricians (see \cite{BM-S}).
The idea is to forget every label and to add a root which is defined
as a star and a branch of this star. We will 
say that a map is rooted at a monomial 
of type $P$ if its root is of type $P$ 
with the marked branch the first one 
in the above construction of a star from a monomial.
We can define for $P$ a monomial, $k_i$ 
a family of integers and $q_i$ a family of
 (this time)  distinct monomials.
\begin{eqnarray*}
\DD_0(P,(q_1,k_1),(q_2,k_2),\ldots,(q_n,k_n))&=\sharp\{&
\textrm{rooted planar  maps with $k_i$ stars of type $q_i$ }\\
& & \textrm{and one of type $P$ which is the root}\, \}
\end{eqnarray*}
and 
$$\DD(P,(q_1,k_1),\ldots,(q_n,k_n))=\sum_{\genfrac{}{}{0pt}{}{1\leq p_i\leq k_i,}{1\leq i\leq n}}\DD(P,(q_1,p_1),(q_1^*,k_1-p_1),\ldots,(q_n,k_n),(q_n^*,k_n-p_n))$$

To go from these rooted maps to the previous one
 we only have to label each star and be careful about
the symmetry of the stars in order to specify a branch by star.
 More precisely, 
let us define  the degree of symmetry $s(q)$
of a monomial $q$  as follows.
Let $\omega:\cxm\ra \cxm$ be the linear function so
that for all $i_k\in\{1,\cdots,m\}$, $1\le k\le p$
$$\omega(X_{i_1}X_{i_2}\cdots X_{i_p})= X_{i_2}\cdots X_{i_p}X_{i_1}$$
and, with $\omega^p=\omega\circ \omega^{p-1}$, define
$$
s(q)=\sharp\{0\leq p\leq \textrm{deg }( q)
-1|\omega^p(q)=q\}.$$ 
We easily see that for all monomial $P$, 
distinct monomials $q_i$ (but eventually, one of
 them may be equal to $P$), and integers $k_i$:

\begin{equation}\label{mapformula}
\DD(P,(q_1,k_1),(q_2,k_2),\ldots,(q_n,k_n))=\frac{\MM(P,(q_1,k_1),(q_2,k_2),\ldots,(q_n,k_n))}{\Pi_{i=1}^nk_i!s_i^{k_i}}
\end{equation}

\subsection{Graphical 
interpretation of Schwinger-Dyson's equations}\label{graphicsec}

We shall now make an assumption on the solutions 
of  Schwinger-Dyson's equation {\bf SD[$V_{\oo{t}}$]}
when the parameters belong to
an open convex neighborhood of the origin,
namely

{\bf (H)} {\it There exists a convex neighborhood $U\in\R^n$, 
 a finite real number $R$ and a family  $\{ \tau_{\oo{t}},
t \in U\}$ of linear forms on $\cxm$ 
so that for all $\oo{t}$ in $U$, $\tau_{\oo{t}}$ is a solution
of {\bf SD[$V_{\oo{t}}$]} which satisfies  {\bf $H( R)$}.}

Note that up to take a smaller set $U$, we can assume 
that the conclusions of Theorem \ref{unique}
are valid, i.e 
  $R\le A(|t|)$ for all $\oo{t}\in U$  since $A(|t|)$
blows up as $|t|$ goes to zero. 

The central result of this article is then

\begin{theo}\label{main}
Assume that (H) is satisfied.
Then
\begin{enumerate}
\item  For any $P\in\cxm$,
$\oo{t}\in U\ra  \tau_{\oo{t}}(P)$ is $\Ca^\infty$ at the origin
in the sense that for all $\oo{k}=(k_1,k_2,\cdots, k_n)\in\N^n$
there exists $\e(k_1+k_2+\cdots+k_n)>0$ 
so that $\partial_{t_1}^{k_1}\cdots\partial_{t_n}^{k_n}
 \tau_{\oo{t}}(P)$
exists on $U_\e=U\cap B(0,\e)$ with $B(0,\e)=
\{\oo{t}\in\R^n:  |t|\le \e\}$.
We let $ \tau^{\oo{k}}(P)=(-1)^{k_1+\cdots+k_n}
\partial_{t_1}^{k_1}\cdots\partial_{t_n}^{k_n}
 \tau_{\oo{t}}(P)|_{\oo{t}=0}$. Then, we have for all $P\in\cxm$
and all $i\in\{1,\cdots,m\}$,
\begin{equation}\label{theequation}
 \tau^{\oo{k}}(X_i P)=\sum_{0\le p_j\le k_j\atop 1\le j\le
n}\prod_{j=1}^n C_{k_j}^{p_j} 
 \tau^{\oo{p}}\otimes  \tau^{\oo{k}-\oo{p}}(D_i P)
+\sum_{1\le j\le n} k_j
 \tau^{\oo{k}-1_j}(\Da_i q_j +\Da_iq_j^*P)
\end{equation}
where $1_j(i)=1_{i=j}$ and  $ \tau^{\oo{k}}(1)=1_{\oo{k}=0}$.
\item For any  monomial 
$P\in\cxm$,   any $k_1,\cdots, k_n\in \N$,
$$ \tau^{\oo{k}}(P)=\MM(P, (q_1,k_1),\ldots,(q_n,k_n)).$$
\end{enumerate}
\end{theo}
\begin{dem}
$\bullet$ The smoothness of $\oo{t}\ra\tau_{\oo{t}}$
comes as in the proof of Theorem \ref{unique}
  from  Schwinger-Dyson's equations and induction on the
degree of the test polynomial function.
Denote $V=V_{\oo{t}}$,
$ \tau= \tau_{\oo{t}}$ and
 take $\oo{t}=(t_1,\cdots,t_n), \oo{t}'=(t_1',t_2',..,t_n')\in U$.
By {\bf SD[V]}, 
$$( \tau_{\oo{t}}- \tau_{\oo{t}' })[(X_i+\Da_i V_{\oo{t}})P]
=(
 \tau_{\oo{t}}-\tau_{\oo{t}'}) \otimes\tau_{\oo{t}}
(D_i P)+\tau_{\oo{t}'}\otimes (
 \tau_{\oo{t}}-\tau_{\oo{t}'})(D_i P)
+  \tau_{\oo{t}' }[(\Da_i V_{\oo{t}'}-\Da_i V_{\oo{t}} )
P] $$
By our finite moment assumption,
we deduce that if $P$ is a monomial function of degree $l-1$,
for any $i\in\{1,\cdots,m\}$,
$$ 
| \tau_{\oo{t}}[(X_i+\Da_i V_{\oo{t}})P]
- \tau_{\oo{t}'}[(X_i+\Da_i V_{\oo{t}})P]|$$
$$
\le 2\sum_{k=0}^{l-2}  
 \max_{Q\mbox{ monomial of degree }\le k}
| \tau_{\oo{t}}[Q]- \tau_{\oo{t}'}[Q]|R^{l-2-k}+  \sum_{1\le i\le n}
|t_i-t'_i| R^{l+D-1}.
 $$
Thus we deduce that for any $p\in\N$,
\begin{eqnarray*}
\D_l( \tau_{\oo{t}},  \tau_{\oo{t}'})&=&
\max_i\max_{P\mbox{ monomial of degree } p-1} |\tau_{\oo{t}}(X_iP)-
\tau_{\oo{t}'}(X_iP)|\\
&\le& 2\sum_{k=0}^{l-2} \D_{k}( \tau_{\oo{t}},  \tau_{\oo{t}'})R^{l-2-k}
+
\sum_{i=1}^n |t_i| \D_{l+d_i -1}( \tau_{\oo{t}},  \tau_{\oo{t}'})
  +\sum_{1\le i\le n}
|t_i-t'_i| R^{ l+  D-1}
.\\
\end{eqnarray*}
Now, let $\gamma\in (0,R^{-1})$
and let's sum both sides of this
inequality multiplied by $\gamma^l$ to obtain, with $d_\gamma( \tau_{\oo{t}},  \tau_{\oo{t}'})
=\sum_{l\ge 0} \gamma^l \D_l( \tau_{\oo{t}},  \tau_{\oo{t}'})$, 
$$d_\gamma( \tau_{\oo{t}},  \tau_{\oo{t}'})
\le 2 (1-\gamma R)^{-1} \gamma^2 d_\gamma( \tau_{\oo{t}},  \tau_{\oo{t}'})
\qquad\qquad\qquad\qquad$$
$$\qquad\qquad+\sum_{i=1}^n|t_i| \gamma^{- d_i+1} d_{\gamma} ( \tau_{\oo{t}},  \tau_{\oo{t}'})
+(1-\gamma R)^{-1} \sum_{1\le i\le n}
|t_i-t'_i| R^{  D-1}
.
$$

Since by definition $\D_l ( \tau_{\oo{t}},  \tau_{\oo{t}'})\le 2R^l$,
$d_\gamma( \tau_{\oo{t}},  \tau_{\oo{t}'})$ is finite for $\gamma R<1$ 
we arrive at 
$$(1-2 \gamma^2(1-R\gamma)^{-1}-  \sum_{1\le i\le n
}|t_i|\gamma^{-D+2}
) d_\gamma( \tau_{\oo{t}},  \tau_{\oo{t}'})
\le (1-R\gamma)^{-1}
\sum_{1\le i\le n}
|t_i-t'_i| R^{  D-1}.$$
Now, for $ |t|$ small enough,
we can find $\gamma=
\gamma( |t|)>0$ so that 
$$1-2\gamma^2 (1-R\gamma)^{-1}- \sum_{1\le i\le n
}|t_i|\gamma^{-D+2}
>0$$
and so 
$$\sum_{l\ge 0} \gamma^l \D_l ( \tau_{\oo{t}},  \tau_{\oo{t}'})
\le C(\oo{t})\sum_{1\le i\le n}
|t_i-t'_i|$$
which implies that for all $l\in\N$
$$\D_l ( \tau_{\oo{t}},  \tau_{\oo{t}'})
\le C(\oo{t})\gamma^{-l} \sum_{1\le i\le n}
|t_i-t'_i|$$
so that for any monomial function $P$,
$\oo{t}\ra   \tau_{\oo{t}}(P)$ is
Lipschitz in $U_\e:=U\cap B(0,\e)$ for $\epsilon$
small enough.
Moreover, we have proved that there exists
$\eta_0(\e)=\gamma^{-1}<\infty$,
 so that
\begin{equation}\label{b1}\D_l ( \tau_{\oo{t}},  \tau_{\oo{t}'})
\le C_0(\e)\eta_0(\e)^l|\oo{t}-\oo{t'}|\mbox{ with }
|\oo{t}-\oo{t'}|=
\max_{1\le i\le n}|t_i-t_i'|.
\end{equation}
Consequently, $ \tau_{\oo{t}}$ is almost surely differentiable
in $U_\e$ and the derivative satisfies
\begin{equation}\label{eq1}
\partial_{t_k} \tau_{\oo{t}}  
[(X_i+\Da_i V_{\oo{t}})P]+ \tau_{\oo{t}}[\Da_i q_k P]
=  \partial_{t_k} \tau_{\oo{t}}\otimes  \tau_{\oo{t}}(D_i P)+\tau_{\oo{t}}\otimes \partial_{t_k} \tau_{\oo{t}}(D_i P)
\end{equation}
for almost all $\oo{t}\in U_\e$.
Since $\cxm$ is countable, these equalities
hold simultaneously for all $P\in\cxm$ almost surely,
let $ U'_\e$ 
be this subset of $U_\e$ of full probability.

(\ref{b1}) implies that
$$\max_{1\le k\le m}\max_{P\mbox{ monomial of degree }l}
|\partial_{t_k} \tau_{\oo{t}}(P)|
\le C_0(\e)\eta_0(\e)^l$$
for 
 all $\oo{t}\in { U'_\e}$. This bound in turn
shows that
we can redo the argument as above to see
that for $|\oo{t}|$ small enough,
$\oo{t}\ra \partial_{t_k}  \tau_{\oo{t}}(P)$
is Lipschitz. Indeed, if we set
$$\D_1(l)=\D_l^1(\partial  \tau_{\oo{t}},\partial  \tau_{\oo{t}}
 )=\max_{1\le k\le m} \max_{P\mbox{ monomial 
of degree } l} |\partial_{t_k} \tau_{\oo{t}}(P)  -
\partial_{t_k} \tau_{\oo{t}'}(P)|$$
we get, for $\oo{t}',\oo{t}\in{ U'_\e}$,
$$\D_1(l)\le 2 \sum_{k=0}^{l-2} \D_1(k) R^{l-2-k} + C_0(\e)
  |\oo{t}-\oo{t'}| l\eta_0(\e)^{l} +
\sum_{i=1}^n |t_i| \D_1(l+d_i -1)$$
so that we get that by summation, for $\gamma<\min(R^{-1}, \eta_0(\e)^{-1})$,
$$(1 -2(1-R\gamma)^{-1}\gamma^2 -\sum_{i=1}^n |t_i|\gamma^{-d_i+1})
 \sum_{l\ge 0}
\D_1(l)\gamma^l \le \gamma^2 C_0(\e)(1-\gamma\eta_0(\e))^{-2}
   |\oo{t}-\oo{t'}|.
$$
Hence, again, we can choose $\eta_1(\e)<\infty$
 big  enough so that
there exists $C_1(\e)<\infty$ so that if $\e$ is small enough
$$\D_1(l)\le C_1(\e)\eta_1(\e)^{l}
 |\oo{t}-\oo{t'}|.$$
In particular, this shows that we can extend 
${\oo{t}}\in { U'_\e} \ra \partial_{t_k}  \tau_{\oo{t}}(P)$
for all monomial functions $P$
continuously in $U_\e$ and so the equality (\ref{eq1}) holds everywhere.
Now, we can proceed by induction to see
that
${\oo{t}}\ra \tau_{\oo{t}}(P)$ 
is $\Ca^{\infty}$ differentiable in a neighborhood of the
origin. More precisely, for any $\oo{k}=(k_1,\cdots, k_n)$
there exists $\e= \e(k_1+k_2+\cdots+k_n) >0$ 
so that on $U_\e$, $$ \tau^{\oo{k}}_{\oo{t}}(P)=
(-1)^{k_1+\cdots +
k_n}\partial_{t_1}^{k_1} \cdots\partial_{t_m}^{k_n}  \tau_{\oo{t}}(P)$$
exists and furthermore satisfies the equation
$$ \tau^{\oo{k}}_{\oo{t}}((X_i+\Da_i V_{\oo{t}})P)=
\sum_{0\le p_i\le k_i\atop 1\le i\le n}\prod_{i=1}^n C^{p_i}_{k_i}
 \tau^{\oo{p}}_{\oo{t}}\otimes  \tau^{\oo{k}-\oo{p}}_{\oo{t}}(D_i P)
+\sum_{1\le j\le m} k_j
 \tau^{\oo{k}-1_j}_{\oo{t}}((\Da_i q_j +\Da_i q_j^*) P)$$
 Applying this result at
the origin, we obtain the announced result.

$\bullet$ We finally show the combinatorial interpretation
of (\ref{theequation}).

Let us first interpret graphically
 $\tau^0=\tau_0$. $\tau_0$ 
satisfies by definition
{\bf SD[0]} which is well known 
to have a unique solution given by the law of $m$ 
free
semi-circular variables (see Voiculescu \cite{Vo0}). 
Then, $ \tau_0(X_{i_1}\cdots X_{i_k})$
can be computed for instance using 
cumulants techniques as developed by R. Speicher \cite{Spe};
it counts the number of planar maps which can be constructed from
the star associated to $X_{i_1}\cdots X_{i_k}$
by gluing together the edges of the star of the same color.
A way to prove  that is to remark 
first that if we have a star
with two branches of the same color,
there is only $1=\tau^0(1)$ ways 
to glue them.
We then proceed by induction over the degree
of the monomial function. We 
let $\Ma(P)$ be the number 
of planar maps with labeled 
stars with a star of type $P$
and  shall show that it satisfies
the same induction relation than $\tau^0(P)$.
Let $i\in\{1,\cdots,m\}$
and $P=X_iQ$. To compute $\Ma(X_iQ)$, we 
 break the edge between the distinguished
branch $X_i$ and the other branch of $Q$ with
which it was glued, then erasing these two branches.
Since the maps are planar, this decomposes the planar map
into two planar maps (see figure \ref{break1})
corresponding respectively 
to the stars $Q_1,Q_2$ for any possible
choices of $Q_1,Q_2$ so that $Q=Q_1 X_i Q_2$.
Hence
$$\Ma(X_iQ)=\sum_{Q=Q_1 X_i Q_2} \Ma(Q_1)\Ma(Q_2).$$
Thus, if $\Ma(R)=\tau_0(R)$ for all 
monomial of degree strictly smaller than $P$,
$$\Ma(X_iQ)=\sum_{Q=Q_1 X_i Q_2} \tau_0(Q_1)\tau_0(Q_2)=
\tau_0\otimes\tau_0(D_iQ)$$
which completes the argument since the right hand side is exactly
$\tau_0(X_iQ)$. 

\begin{figure}[h!]
\psfrag{x1}{\begin{LARGE}{$X_{1}$}\end{LARGE}}
\psfrag{x2}{\begin{LARGE}{$X_{2}$}\end{LARGE}}
\psfrag{p1}{\begin{LARGE}{$X_{1}X_2^2X_1$}\end{LARGE}}
\psfrag{p2}{\begin{LARGE}{$X_1^2X_{2}^2$}\end{LARGE}}
\psfrag{m1}{\begin{LARGE}{Edge to be broken}\end{LARGE}}
\psfrag{m2}{\begin{LARGE}{Boundary of the maps}\end{LARGE}}
\begin{center} \resizebox{!}{8cm}{\includegraphics{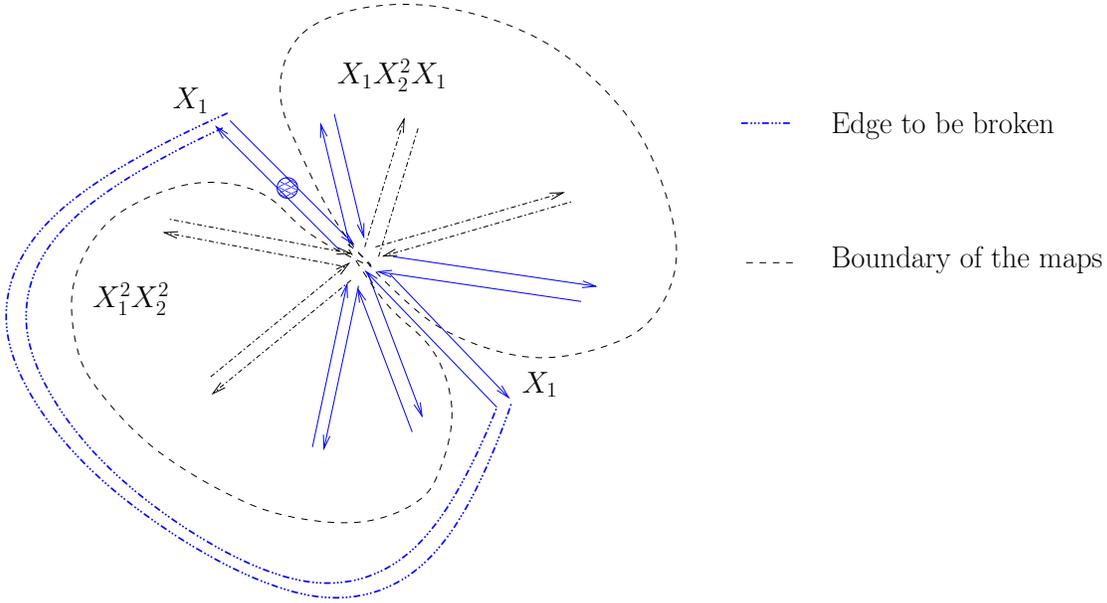}}
\end{center}
\caption{The decomposition  $P(X)=X_1X_2^2X_1^4X_2^2$ into $X_1X_2^2X_1
\otimes X_1^2X_2^2$ }\label{break1}
\end{figure}

We now consider the general case; let us
assume   that 
for $|\oo{k}|\le M$, the graphical interpretation has been obtained for
all monomial and that for $|\oo{k}|= M+1$, it has been proved for monomial
of degree smaller are equal to $L$. By the preceding, we can take $M\ge 1$
and $L\ge 1$ since  for all $\oo{k}\neq 0$, $ \tau^{\oo{k}}(1)=0$.
Again, we shall
show that ${\mathcal M}( P, (q_1,k_1),\cdots,(q_n,k_n))$
satisfies the same induction relation than $\tau^{\oo{k}}(P)$.

Let us consider a star of type $X_i P$ (rooted at the branch $X_i$, with its inner orientation)
with $P$ a monomial of degree less than $L$ and $|\oo{k}|=\sum k_i=M+1$. 
Now, in order to compute ${\mathcal M}( X_iP, (q_1,k_1),\cdots,(q_n,k_n))$,
we break  the edge between
the distinguished branch
$X_i$  (which has color $i$)
and the other branch with which 
it was glued. 

The first possibility is that it was glued with
an edge of the star $P$. Then, since the
maps are planar, this decomposes the map
in  two planar maps. If this branch was given by the $X_i$ 
so 
that $P=P_1X_iP_2$,  one of this planar map contain the star of type $P_1$
and the other the star of type $P_2$, which have also a distinguished branch  and
are  oriented.
If one of this planar map
is glued with $ k_j$ stars
of type $q_j$ or $q_j^*$ , $0\le k_j\le n$,
the other map is glued with the remaining
stars, that is $k_j-p_j$ 
stars of type
$q_i$ or $q_i^*$. There are $\prod_{j=1}^n C_{k_j}^{p_j}$
ways to choose $p_j$ among $k_j$ 
stars of type $q_j$ or $q_j^*$ for $1\le j\le n$
(recall here that stars are labeled).
Since we do that for all $(P_1,P_2)$ so that $P$ have the above 
decomposition,
we obtain the planar  maps corresponding actually 
to the stars associated with the monomials of  $D_iP$.
Note that the case where one of the monomial
in $D_iP$ is the monomial
$1$ shows up  when
$P=X_iQ$ or $Q X_i$ for some monomial $Q$
and the weight corresponds
then to the case where we glue 
the first branch $X_i$ in $X_iP$
with its left or right neighbor.
In this case, none of these two branches 
can be glued with another star,
and there is only one possibility
to glue these two branches otherwise,
which corresponds to the weight
$ \tau^{\oo{k}}(1)=1_{\oo{k}=0}$.

Hence, the number of planar maps corresponding to this
configuration 
 is given by 

$$\sum_{0\le p_j\le k_j\atop 1\le j\le n}\sum_{P=P_1X_iP_2} \prod_{1\le j\le n}
C_{k_j}^{p_j}\Ma(P_1, (q_1,p_1),\cdots, (q_n,p_n))\Ma(P_2, (q_1,k_1-p_1),\cdots, (q_n,k_n-p_n))$$
$$=\sum_{0\le p_j\le k_j\atop 1\le j\le n} \prod_{1\le j\le n}
C_{k_j}^{p_j}
 \tau^{\oo{p}}\otimes  \tau^{\oo{k}-\oo{p}}(D_i P)$$
where we finally used our induction hypothesis.

The other possibility is that this edge is
glued with a star of type $q_j^{\epsilon}$ for some $j\in\{1,\cdots, n\},\epsilon\in\{.,*\}$.
In this case, erasing the edge means that
we destroy a star of type $q_j^{\epsilon}$
and replace the stars of type $X_iP$ and $q$
glued together with a single star $P$ glued with the star
$q_j^{\epsilon}$ in place of $X_i$ with an edge of color $i$ removed;
if $q_j^{\epsilon}=Q_1X_i Q_2$, we replace  the two stars of type
 $X_iP$ and $q_j^{\epsilon}$ by a single one of type $Q_2Q_1 P$
(see figure \ref{break2}).
Since we do that with all the possible edges of color $i$ in $q_j^{\epsilon}$,
we find that we can glue 
all monomials appearing in $\Da_i q_j^{\epsilon}$,
and so the corresponding weight is given by
$ \tau^{\oo{k}-1_j}(\Da_i q_j^{\epsilon} P)$ times $k_j$, the number
of ways to choose one star among $k_j$
of type $q_j^{\epsilon}$.

\begin{figure}[ht!]
\psfrag{x1}{\begin{LARGE}{$X_{1}$}\end{LARGE}}
\psfrag{x2}{\begin{LARGE}{$X_{2}$}\end{LARGE}}
\begin{center} \resizebox{!}{7cm}{\includegraphics{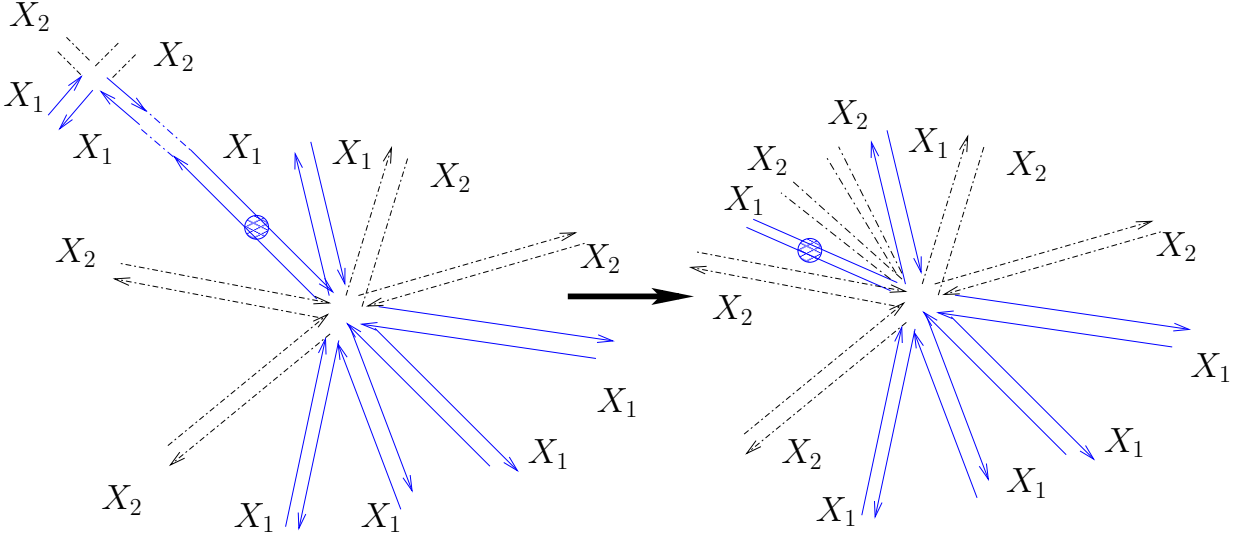}}
\end{center}
\caption{The merging  of $q(X)=X_1^2X_2^2=X_1 X_1 X^2_2$ and 
$X_1P$ into $ X_1 X^2_2 P$
 }\label{break2}
\end{figure}

Hence, by induction, we proved
that the number of planar
maps with $k_j$ stars of type $q_j$ or $q_j^*$
and one of type $X_iP$ is given 
by

\begin{align}
\Ma((X_iP,1), (q_1,k_1),..., (q_n,k_n)) &=\sum_{0\le p_j\le k_j\atop 1\le j\le m}\prod_{j=1}^n C_{k_j}^{p_j}
 \tau^{\oo{p}}\otimes  \tau^{\oo{k}-\oo{p}}(D_i P)
+\sum_{1\le j\le m} k_j
 \tau^{\oo{k}-1_j}(\Da_i (q_j +q_j^*)P)\\
&= \tau^{\oo{k}}(X_i P)
\end{align}
for all $i\in\{1,\cdots,m\}$.
This shows that the graphical interpretation holds for
all $L$ and $|\oo{k}|\le M+1$.
 We can 
start the induction
since we know that $ \tau^{\oo{k}}(1)=1_{\oo{k}=0}$.
This completes the proof.

\end{dem}

\noindent
{\bf Remarks:}
\begin{enumerate}
\item 
This graphical interpretation can be sometimes simplified for particular $V$.
For example, consider $V_{t,u,c}=tA^4+uB^4+cAB$ which appears in the Ising model. First, one may notice 
that this is not in the form of the theorem but we can replace 
$\tr (V)$ by
$$\tr(\frac{t}{2}(A^4+A^4)+\frac{u}{2}(B^4+B^4)+\frac{c}{2}(AB+BA))$$
so that the theorem can be applied.
Now, one may see that there will some redundancy in the enumeration of maps
given by this potential as for example vertices of type $AB$ are isomorphic to
vertices of type $BA$.
But there will be some simplification with the factors $\frac{1}{2}$ so that 
finally the theorem will give, for all monomial $P$, if $\tau_{t,u,c}$ is the a solution to  {$\bf SD[V_{t,u,c}]$} given by the theorem then

$$(-1)^{m+n+r}
\partial_{t}^{m}\partial_{u}^{n}\partial_{c}^{r}
 \tau_{t,u,c}(P)|_{t=u=c=0}=\MM_0(P,(A^4,m),(B^4,n),(AB,r))$$

\item 
Note that this 
graphical approach can be generalized to matrix models 
with more complex
 potentials involving tensor 
products. For example,
one  can consider a potential $V$ which is a sum of monomials and
of tensor products of monomials:
$$V_{\oo{t}}=\sum_i t_i (q_i^1\otimes \cdots \otimes q_i^{d}+
(q_i^1)^*\otimes \cdots \otimes (q_i^{d})^*)$$
and the associated measure with
density with respect to $\mu_N^{\ot m}$ given by
$Z_N^{-1}e^{-N^{2-d}(\tr)^{\ot d} V_{\oo{t}}}$.
Then one can write the generalized Swinger Dyson's equation:

$$\tau\otimes \tau(D_i P)=\tau( X_iP)+\sum_{k,j}t_k
 \tau^{\otimes d_k}(q_k^1\otimes\cdots
\otimes \Da_iq_k^{j} P\otimes\cdots\otimes$$
$$q_k^{d}+
(q_k^1)^*\otimes\cdots
\otimes \Da_i(q_k^{j})^* P\otimes\cdots\otimes (q_k^{d})^*)$$

The previous results remain
 valid up to a graphical interpretation 
 of the new term.
For example 
$q^1\otimes\cdots\otimes q^k$ will be a bunch of $k$ loops, the first one containing 
the branches of $q^1$ in the clockwise order, the second one the branches of $q^2$ ...
Note that one of these $k$ loops is redundant, so that we can choose to delete one loop
in the graphical representations.
These vertices will split the map into parts as the vertices which will be placed in a loop
can not be linked to any vertices in an other loop.
\end{enumerate}

\subsection{Existence of an analytic solution to 
 Schwinger-Dyson's equation}\label{existanal}
The aim of this section is to prove that for all monomials 
$(q_j)_{1\le j\le n}$,
there exists a convex neighborhood of the origin
(actually an open ball) and a finite constant $R$  so that 
 hypothesis {\bf (H)}  of section \ref{graphicsec}
is satisfied.
Moreover, we show it depends analytically on ${\oo{t}}$ in
a neighborhood of the origin.
Let $V_{\oo{t}}$ be as before. 
\begin{theo}\label{existanaltheo}
There exists an open  neighborhood $U\subset{\mathbb R}^n$ 
of the origin
(a ball of positive radius)
  such that for ${\oo{t}}\in U$, there exists  $\tau_{\oo{t}}\in\cxm^*$ 
satisfying {\bf SD[$V_{\oo{t}}$]} such that:
\begin{itemize}
\item
${\oo{t}}\ra \tau_{\oo{t}}$ is analytic on $U$, i.e. there exists
 $ \tau^{\oo{k}},{\oo{k}}\in{\NN}^n$ in $\cxm^*$ 
such that for all $P$ in $\cxm$, $t$ in $U$,
$$\tau_{\oo{t}}(P)=\sum_{{\oo{k}}\in{\NN}^n}\prod_{1\le i\le n}
\frac{(-t_i)^{k_i}}{k_i!}
 \tau^{\oo{k}}(P)$$
and the serie converges absolutely on $U$.
\item
$ \tau^{\oo{k}}(P)=(-1)^{\Sigma k_i}\partial_{t_1}^{k_1}\cdots\partial_{t_n}^{k_n}
\tau_{\oo{t}}(P)|_{\oo{t}=0}=\MM(P,(q_1,k_1),\cdots,(q_n,k_n))$

\item There exists $R<\infty$ so that for all ${\oo{t}}\in U$,
all $i_i\cdots i_l\in \{1,\cdots,m\}^l$,
all $l\in\N$,
$$|\tau_{\oo{t}}(X_{i_1}\cdots X_{i_l})|\le R^l.$$
\end{itemize}
\end{theo}

\nn
{\bf Remark: }
It will be useful to use sometimes for $\tau_{\oo{t}}$  an alternative expression
related to the enumeration of rooted maps.
Using (\ref{mapformula}), one can obtain inside the domain of convergence, for all monomial $P$:
$$\tau_{\oo{t}}(P)=\sum_{{\oo{k}}\in{\NN}^n}\prod_{1\le i\le n}
(-s(q_i)t_i)^{k_i}\DD(P,(q_1,k_1),\cdots,(q_n,k_n)).$$
\begin{dem}
If we have such a solution, it satisfies  assumption
 {\bf (H)}  and by $(\ref{theequation})$
if $\oo{k}!=\prod k_i!$,
$$\frac{ \tau^{{\oo{k}}}(X_iP)}{\oo{k}!}=
\sum_{0\le p_j\le k_j\atop 1\le j\le n}\sum_{ P=P_1X_iP_2}
\frac{ \tau^{\oo{p}}(P_1)}{\oo{p}!}\frac{\tau^{\oo{k}-\oo{p}}(P_2)}{(\oo{k-p})!}
+\sum_{1\le j\le m \atop k_j\neq 0} 
\frac{ \tau^{\oo{k}-1_j}(\Da_i (q_j + q_j^*)P)}{(\oo{k}-1_j)!}$$
where the second sum runs over all monomials $P_1,P_2$ so that
$P$ decomposes into $P_1X_iP_2$.
We can use this formula to define the $\tau^{\oo{k}}$ by induction, the graphical interpretation is directly satisfied.
\par

We must control the growth of the $\tau^{\oo{k}}(P)$'s.
Our induction hypothesis will be
that for $\oo{k}$ so that $\sum_i k_i\le M-1$
and all monomial $P$, as well
as for $\sum k_i=M$ and monomials $P$  of degree smaller than $L$,
$$\left|\frac{\tau^{\oo{k}}(P)}{\oo{k}!}\right|\leq
 A^{\sum k_i}B^{deg P}\prod_i C_{k_i}C_{deg P}$$
where the $C_k$ are the Catalan's numbers which satisfy
\begin{equation}\label{catalan}
C_{k+1}=\sum_{p=0}^k C_pC_{k-p},\quad
C_0=1,\quad \frac{C_{k+l}}{C_l}\le 4^k \quad \forall l,k\in\N.
\end{equation}
Here,  $\deg P$ denotes the degree of the
monomial $P$ and we can assume $B\ge 2$ without loss of
generality.
Our induction is trivially true for $\oo{k}=0$ and all $L$ since
$\tau^{\oo{0}}=\s^m$ is the law of $m$ free semi-circular
variables which are uniformly bounded by $2$ so that
$$|\tau^{\oo{0}}(P)|\le 2^{deg P}$$
Moreover, it is satisfied for
 all $\oo{k}$ and
$L=0$ since then
$ \tau^{\oo{k}}(1)=1_{\oo{k}=0}$.
Let us assume that it is true for all $\oo{k}$ such that
$\sum k_i\le M-1$ and all monomials, and for $\oo{k}$
such that $\sum k_i=M$ and monomials $P$ of degree less 
than $L$ for some $L\ge 0$. Then

\begin{eqnarray*}
\left|\frac{ \tau^{\oo{k}}(X_iP)}{k!}\right|&\leq &
\sum_{0\le p_j\le k_j\atop 1\le j\le n}\sum_{ P=P_1X_iP_2}
A^{\sum k_i}B^{deg P-1}\prod_{i=1}^n C_{p_j}C_{k_j-p_j}C_{deg P_1}C_{deg P_2}\\
& &+2\sum_{1\le l\le n}A^{\sum k_j-1}\prod_j C_{k_j}B^{deg P+deq q_l-1}C_{deg P+deq q_l-1}\\
&\leq& A^{\sum k_i}B^{deg P+1}\prod_i C_{k_i}C_{deg P+1}
\left(\frac{4^{n}}{B^2}+2\frac{\sum_{1\le j\le n}
B^{deg q_j-2}4^{deg q_j-2}}{A}\right)
\end{eqnarray*}
where we used (\ref{catalan}) in the last line.
It is now sufficient to choose $A$ and $B$ such that 
$$\frac{4^{n}}{B^2}+2\frac{\sum_{1\le j\le n}B^{deg q_j-2}
4^{deg q_j-2}}{A}\le1$$
(for instance $B=2^{n+1}$ and $A=4nB^{D-2}4^{D-2}$)
to verify  the  induction hypothesis
works for polynomials of all degrees (all $L$'s). 
\par

Then
$$\tau_{\oo{t}}(P)=\sum_{k\in{\NN}^n}\prod\frac{(-t_i)^{k_i}}{k_i!}\tau^{\oo{k}}(P)$$
is well defined for $|t|<(4A)^{-1}$. Moreover, for all  monomial $P$,
$$|\tau_{\oo{t}}(P)|\le \sum_{k\in\N^n} \prod_{i=1}^n (4t_i A)^{k_i} (4B)^{deg P}
\le \prod_{i=1}^n (1- 4At_i)^{-1} (4B)^{deg P}.$$
so that for small $t$, $\tau_{\oo{t}}$ has an uniformly bounded support.

\end{dem}

Hence, we see that the enumeration of planar maps could be reduced to the
study of  Schwinger-Dyson's equations {\bf SD[V]}.
For instance, the asymptotics of such enumeration
can be obtained by studying the optimal domain in which 
the solutions are analytic.
Matrix models can be useful to study also
the solution.
For instance, we shall deduce from this
approach that the solutions to {\bf SD[V]} are
tracial states (the positivity condition 
being unclear a priori).

\section{Existence of tracial states solutions to 
 Schwinger-Dyson's equations from matrix
models}\label{mm}
We let $V=V_{\oo{t}}$ be a polynomial 
function as before,
and consider
$$Z^N_V=\int e^{-N\tr(V(A_1,\cdots, A_m))}
\mu_N(dA_1)\cdots \mu_N(dA_m)$$
and $\mu^N_V$  the associated Gibbs
measure
$$\mu^N_V(dA_1,\cdots,dA_m)= (Z^N_V)^{-1} e^{-N\tr(V(A_1,\cdots, A_m))}
\mu_N(dA_1)\cdots \mu_N(dA_m).$$
This is well defined provided 
that we assume that the monomials 
of highest degree in $V_{\oo{t}}$ are sufficiently
large to make $Z^N_V$
finite. We shall assume for instance  that
\begin{equation}\label{ass}
V_{\oo{t}}(\bX)=\sum_{1\le i\le n} t_i (q_i(\bX)+q_i^*(\bX))+
\sum_{n+1\le i\le n+m} t_i X_{i-n}^{D} \end{equation}
with $D$ even and monomial functions $q_i$ of degree less or equal than
$D-1$ and $t_i>0$ for $i\in\{n+1,\cdots,n+m\}$.
We shall see in the last paragraph of
this section that such assumption can be removed
provided a cut-off is added.

The empirical distribution of $m$ matrices $\bA=(A_1,\cdots,A_m)\in\Ha_N^m$
is defined
as the element of $\Ma^m_{ST}$ such that
$$\mun_{A_1,\cdots,A_m}(F)=\mun_{\bM}(F)=\frac 1N\tr(F( A_1,\cdots,A_m))$$
for all $F\in\Ca_{st}^m(\C)$.
Note that the empirical distribution could be defined 
as well as an element 
of $\Ma^m$ but since the random matrices $(A_1,\cdots,A_m)$
under $\mu^N_V$ have a priori no uniformly bounded
spectral radius, the topology of weak convergence would not
be suitable then.

We shall see that if we know that
a limit point of $\mun_{\bA}$ under
$\mu^N_V$ are compactly supported,
then it satisfies {\bf{SD}[V]}.
In a second part, we shall give examples
of potential $V$ for which this assumption is
satisfied. Finally, we discuss localized matrix
integrals and show that bounded solutions to  {\bf SD[V]}
for small potentials can always be constructed by
localized matrix integrals.

\subsection{Limit points of empirical distribution of matrices
following matrix models satisfy the {\bf SD[V]} equations}\label{limitpoint}

We claim that

\begin{theo}\label{mmSD}
Assume \eqref{ass}. Then 
\begin{enumerate}
\item
There exists $M<\infty$ so that
$$\limsup_{N\ra\infty}\mun_{A_1,\cdots, A_m}
(X_i^{D})\le M$$
$\mu^N_V$ almost surely for all $i\in\{1,\cdots,m\}$.
\item The limit points of $\mun_{A_1,\cdots, A_m}$
for the $\Ca_{st}^m(\C)$-topology 
satisfy the `weak' Schwinger-Dyson equation
$$\tau\otimes \tau(D_i F)=\tau( (\Da_i V+X_i)F)\quad {\bf(WSD)[V]}$$
for all $F\in\Ca_{st}^m(\C)$.
\end{enumerate}
\end{theo}
Note here that $ (\Da_i V+X_i)F$ does not belong to
$\Ca_{st}^m(\C)$ so that it is not clear what {\bf (WSD)[V]}
means a priori. We define it by the following;
there exists a sequence $V^\d\in \Ca_{st}^m(\C)$
so that
$$\lim_{\d\ra 0}
\max_{1\le i\le m}\sup_{\tau(X^{D}_i)\le M}\tau ( |\Da_i V^\d-\Da_i V-X_i|)=0$$
from which, since any $F\in  \Ca_{st}^m(\C)$
is uniformly bounded, 
$$\lim_{\d\ra 0}
\max_{1\le i\le m}\sup_{\tau(X^{D}_i)\le M}|\tau(F \Da_i V^\d)-\tau(
 F( \Da_i V+X_i))|=0$$
is well defined. 

\nn
\begin{dem}
$\bullet$ The first point is trivial
since by Jensen's inequality,
$$Z_N^V\ge \exp\{-N^2 \int \frac{1}{N}\tr(V(\bA))\prod_{1\le i\le m} 
d\mu_N(A_i)\}\ge \exp\{cN^2\}$$
for some $c>-\infty$, 
where the last inequality comes from the fact that (see \cite{Vo0})
$$\lim_{N\ra\infty}\int \frac{1}{N}\tr(V(\bA))\prod_{1\le i\le m} 
d\mu_N(A_i)=\s^m( V)<\infty$$
where $\s^m$ is the law of $m$ free semi-circular variables.

Now, observe that by H\"older's inequality,
$$ |\mun_{\bA}(q_i)|\le \max_{1\le i\le m}
\mun_{\bA}(|X_i|^{D-1}+1)$$
so that we deduce 
$$\mun_\bA(V)\ge \sum_{i=1}^m \left( t_{i+n} \mun_{\bA}(A_i^{D}
) -c(\oo{t}) \mun_{\bA}(|A_i|^{D-1})-c(\oo{t})\right)$$
with a finite constant $c(\oo{t})$.
Since $t_{i+n}>0$, we conclude that 
$\mun_\bA(V)\ge m|t|M/2$ when 

\noindent
 $\max_{1\le i\le m}  
\mun_{\bA}(A_i^{D}
)\ge M$ for $M$ large enough. 
Thus 
\begin{equation}\label{expestimate}
\mu^N_V\left(\max_{1\le i\le m}  
\mun_{\bA}(A_i^{D}
)\ge M\right) \le 
 e^{-2^{-1} N^2 Mm|t|}e^{-cN^2}
\end{equation}
goes to zero exponentially fast when $M>\frac{2c}{m|t|}$. The claim follows by
Borel-Cantelli's lemma.

$\bullet$ We proceed as in \cite{CDG2}, following 
a common idea in physics, which is to make, in $Z^N_V$,  the
change of variables $X_i\ra X_i +N^{-1} F(\bX)$ 
for a given $i\in\{1,\cdots,m\}$ and $F\in \Ca_{st}^m(\R)$.
Noticing that the Jacobian for this change of variable is
$$|J|= e^{ N\mun_{\bA}\otimes \mun_{\bA}( D_i F) +O(1)}$$
we get
that
$$\int e^{ (N\mun_{\bA}\otimes \mun_{\bA}( D_i F)
-N^2\mun_{\bA} ( N^{-1} X_i F(\bX)+V( X_i +N^{-1} F(\bX))- V(X_i))
} \mu^N_V(d\bA)=O(1)$$

from which we deduce that
$$\int_{\max_{1\le i\le m}  
\mun_{\bA}(A_i^{D}
)\le M
}e^{ (N\mun_{\bA}\otimes \mun_{\bA}( D_i F)
-N^2\mun_{\bA} ( N^{-1} X_i F(\bX)+V( X_i +N^{-1} F(\bX))- V(X_i))
} \mu^N_V(d\bA)=O(1).$$
Hence, we conclude by Chebychev inequality and (\ref{expestimate})
that  for $M$ big enough, any $\d>0$, there exists $\eta>0$,
so that if we denote
$$E_N=\mun_{\bA}\otimes \mun_{\bA}( D_i F)
-N\mun_{\bA} ( N^{-1} X_i F(\bX)+V( X_i +N^{-1} F(\bX))- V(X_i))$$
then
$$\mu_V^N\left( \{\max_{1\le i\le m}  
\mun_{\bA}(X_i^{D}
)\le M\}\cap \{|E_N|\leq\d\}\right)\ge 1-e^{-\eta N}.$$

Moreover,
$$\mun_{\bA}(V( X_i +N^{-1} F(\bX))- V(X_i)))
= N^{-1} \mun_{\bA}(\Da_i V F)+ R_N$$
with a rest $R_N$ of order $N^{-2} \max_{1\le i\le m}  
\mun_{\bA}(X_i^{D-2})$ which we can neglect 
on \hfill

\nn
$\max_{1\le i\le m}  
\mun_{\bA}(A_i^{D})\le M$.
This shows, by Borel-Cantelli's Lemma, that for all $F\in \Ca_{st}^m(\R)$, 
$$\mun_{\bA}\otimes \mun_{\bA}( D_i F)-
\mun_{\bA} (X_i F+\Da_i V F)$$
goes to zero almost surely. This result extends to $F\in
\Ca_{st}^m(\C)$ since it can always be decomposed into
the sum of two elements of $\Ca_{st}^m(\R)$.
Moreover, 
 if we let $A_i^\epsilon=A_i (1+\e A_i^2)^{-1}= A_i (\sqrt{-1}
+\sqrt{\e} A_i)^{-1}(-\sqrt{-1}
+\sqrt{\e} A_i)^{-1}\in \Ca_{st}^m(\C)$, 
then
again by H\"older's inequality
$\tau
(|\Da_i V(A_i)-\Da_i V(A_i^\e)|)$  goes to zero uniformly 
on $\max_{1\le i\le m}  
\tau(A_i^{D}
)\le M$. This shows that $\mu\ra\mu( (\Da_i V+X_i)F)$
is continuous for the weak $\Ca_{st}^m(\C)$-topology on $\{\mu( A_i^{D})\le M\}$
for any $F\in\Ca_{st}^m(\C)$.
 Therefore, since 
$\Ma^m_{ST}$ is compact, we
conclude that any limit point of $\mun_{\bA}$ 
satisfies 
$$\tau\otimes \tau( D_i F)=
\tau ((X_i +\Da_i V )F)$$

\end{dem}

We therefore have
the
\begin{cor}\label{comp}
Assume that there exists a limit point $\tau_{V}$ 
of $\mun_{\bA}$ under $\mu^N_V$ which is
compactly supported.
Then, it satisfies  Schwinger-Dyson's equation
{\bf SD[V]}.
\end{cor}
\begin{dem} 
The proof is straightforward since 
if $\tau_{V}$ is compactly supported
it is equivalent to say that $\tau_V$ satisfies 
 {\bf WSD[V]} or {\bf SD[V]} since 
$\Ca_{st}^m (\C)$ is dense in the set of polynomial functions
(approximate the $A_i$'s by the $A_i^\d$'s
defined in the previous proof).
\end{dem}
Let us also give the final
argument to deduce
convergence of the free energy
from the previous considerations.

\begin{theo}\label{convfe}
\begin{enumerate}
\item
Assume that $\mun_{\bA}$
converges in $\Ma^m_{ST}$ almost surely 
or in expectation under $\mu^N_{V_{\oo{t}}}$
towards $\tau_{\oo{t}}$  solution to {\bf SD[$V_{\oo{t}}$]}
for $\oo{t}$ in 
a convex  neighborhood $U$  of the
origin. Assume
furthermore that $\max_p\mu^N_{V_{t_1,\cdots,t_{k-1}, s, 
t_{k+1},\cdots, t_n}} (\mun_{\bA}(|X_p|^l))$
is uniformly bounded for 

\nn
$ 
(t_1,\cdots,t_{k-1}, s, t_{k+1},\cdots, t_n)\in U$
and $N$ large enough for some $l$ strictly greater than the
degree of $V_{\oo{t}}$.
 Then, $$
F_{V_{\oo{t}},k}^N=N^{-2}
\log\left( \frac{Z^N_{V_{\oo{t}}}}{Z^N_{V_{(t_1,
\cdots, t_{k-1}, 0,t_{k+1}..t_n)}}}\right)$$
converges as $N$ goes to infinity
towards a limit $F_{V_{\oo{t}},k}$.
Moreover,
$$F_{V_{\oo{t}},k}=-\int_0^{t_k}
\tau_{(t_1,\cdots,t_{k-1},s,t_{k+1},\cdots,t_n)}(q_k+q_k^*) ds.$$

If furthermore $\tau_{\oo{t}}$ is uniformly compactly supported 
in $U$, we deduce that
$\oo{t}\ra F_{V_{\oo{t}}}$ is $\Ca^\infty$ 
in a neighborhood of the origin
and $(-1)^{\sum p_i}\partial_{t_1}^{p_1}\cdots \partial_{t_n}^{p_n} F_{V_{\oo{t}},k}|_{\oo{t}=0}$
is the number of planar maps  with $p_i$ stars of type $q_i$ or $q_i^*$ 
when $p_k\ge 1$.

\item Assume that for $\oo{t}$ in a open convex neighborhood $U$
of the origin the limit points of $\mun_{\bA}$
under $\mu^N_{V_{\oo{t}}}$ are uniformly compactly supported.
Assume further
that $\max_p|\mu^N_{V_{\oo{t}}} (\mun_{\bA}(|X_p|^l))|$
is uniformly bounded (independently of $\oo{t}\in U$) for $N$
 large enough and some $l$ strictly larger than the degree 
of $V_{\oo{t}}$.
Then, $\mun_{\bA}$ converges $\mu^N_{V_{\oo{t}}}$-almost surely towards 
$\tau_{\oo{t}}$ described  in Theorem \ref{existanaltheo}
for $\oo{t}\in U\cap B(0,\e)$ for some $\e>0$ small enough
and for $\oo{t}\in U\cap B(0,\e)$, 
$$F_{V_{\oo{t}}}^N=N^{-2}
\log( {Z^N_{V_{\oo{t}}}})$$
converges as $N$ goes to infinity towards
$$F_{V_{\oo{t}}}=\sum_{\oo{k}\in \N^n\backslash (0,..,0)}
\prod_{1\le i\le n} \frac{(-t_i)^{k_i}}{k_i!}
\Ma((q_1,k_1),\cdots,(q_n,k_n)).$$

\end{enumerate}
\end{theo}

Note above that the last serie as
a positive radius of convergence according to
Theorems \ref{main} and \ref{existanaltheo}. This emphasizes that
the possible divergence of $F^N_{V_{\oo{t}}}$ does not survive the 
large $N$ limit.

\noindent
\begin{dem}
$\bullet$
By differentiating $N^{-2}\log Z^N_{V_{\oo{t}}}$
with respect to $t_k$ we obtain that
$$\partial_{t_k}N^{-2}\log
Z^N_{V_{\oo{t}}}=-\mu^N_{V_{\oo{t}}}(\mun_\bA(q_k+q_k^*)).$$
But,  under assumption, $(\mun_\bA(q_k+q_k^*))_{N\in\N}$
converges almost surely and is uniformly integrable
so that $\mu^N_{V_{\oo{t}}}(\mun_\bA(q_k+q_k^*))$
is a uniformly bounded sequence which 
converges as $N$ goes to infinity towards $\tau_{\oo{t}}(q_k+q_k^*)$
for $\oo{t}\in U$. Integrating with respect to $t_k$ yields
the convergence with $F_{V_{\oo{t}}}$ as above
by dominated convergence theorem.
The last part of the first point theorem is a direct consequence of
Theorem \ref{main}.

$\bullet$ By Corollary \ref{comp} and Theorem \ref{unique},
we see that our hypothesis implies that for $\oo{t}\in U\cap B(0,\e)$
for some $\e>0$, the limit points 
of $\mun_\bA$ are unique and given by $\tau_{\oo{t}}$.
Hence, $\mun_\bA$ converges in $\Ma^m_{ST}$ almost surely
towards $\tau_{\oo{t}}$. Since we assumed
our family uniformly integrable, we deduce
that $\mu^N_{V_{\oo{t}}}(\mun_\bA(q_l+q_l^*))$
converges as $N$ goes to infinity
towards $\tau_{\oo{t}}(q_l+q_l^*)$ for all $l\in\{1,\cdots, n\}$
and we see as above that for all $i\in \{1,\cdots, n-1\}$,
$$\frac{1}{N^2}
\log\left( \frac{Z^N_{V_{(0,..,0,t_i..t_n)}}}{Z^N_{V_{(0,..,0,t_{i+1}..t_n)}
}}\right)$$
converges as $N$ goes to infinity
towards a limit $$F_{V_{\oo{t}},i}=-\int_0^{t_i}
\tau_{(0,\cdots,0,s,t_{i+1},\cdots,t_n)}(q_i+q_i^*) ds.$$
Hence,
since we know that $F_{V_0}=1$,
\begin{eqnarray*}
\lim_{N\ra\infty}
\frac{1}{N^2}\log Z^N_{V_{\oo{t}}}&=&\sum_{i=1}^{n}
\lim_{N\ra\infty}\frac{1}{N^2}\log 
\left( \frac{Z^N_{V_{(0,..,0,t_i..t_n)}}}{Z^N_{V_{(0,..,0,t_{i+1}..t_n)}
}}\right)\\
&=&
-
\sum_{i=1}^{n}\int_0^{t_i}
\tau_{(0,\cdots,0,s,t_{i+1},\cdots,t_n)}(q_i+q_i^*) ds
\\
&=& 
\sum_{i=1}^{n}\sum_{k_i,..,k_n\in \N^{n-i}}
\prod_{i+1\le j\le n} \frac{(-t_j)^{k_j}}{k_j!}
\frac{(-t_i)^{k_i+1}}{(k_i+1)!}\tau^{(0,..,0,k_i,\cdots,k_n)}
(q_i+q_i^*)\\
\end{eqnarray*}
where we used in the last line 
Theorem \ref{existanaltheo}.
Noting that $\tau^{(0,..,0,k_i,\cdots,k_n)}
(q_i+q_i^*)$ is the number of planar maps with
$k_j$ stars of type $q_j$ or $q_j^*$ for $j\ge i+1$
and $k_i+1$ stars of type $q_i$ or $q_i^*$ , we conclude the 
proof.

\end{dem}

We shall in the next section provide 
a generic example where the 
assumption of the second point of Theorem \ref{convfe}
is satisfied (in fact, a slightly different version
since we do not prove that the almost sure
limit points of $\mun_\bA$ satisfy our
compactness assumption, but their average
do, which still guarantees the result).

\subsection{Convex interaction models}\label{convexsec}
Let us assume that we consider 
a  matrix model with potential $V$
such that 
\begin{equation}\label{conv}
\phi_{V,a}^N: (A_k({ij}))\in (\R^{N^2})^m
\ra \tr( V(A_1,\cdots,A_m))- \frac{a}{2}\sum_{k=1}^m 
\tr (A_k^2)
\end{equation}
is convex in all dimensions for some $a<1$,
i.e the Hessian of $\phi_{V,a}^N$ 
is non negative for all $N\in\N$.
An example is 
$V$  of the form
$$V(A_1,\cdots,A_m)=\sum_{i=1}^n t_i (\sum_{k=1}^m
\alpha_k^i A_k)^{2p_i}$$
with non-negative  $t_i$'s, integers $p_i$'s
and real $\alpha$'s.  Indeed, by Klein's lemma (c.f.
\cite{GZ}), since $x\ra (\sum\alpha_k x_k)^{2p_i}$
is convex, 
$$\bA \ra\tr (\sum\alpha_i A_i)^{2p_i}$$
is also convex (Here $\bA$,
by an abuse of notations, denotes the entries of the $m$-uple of matrices
$\bA=(A_1,\cdots,A_m)$). 

Then, we shall prove that
\begin{theo}\label{convex}
Let $V$ be a  self-adjoint polynomial
function which satisfies (\ref{conv}). Then
\begin{itemize}

\item 
There exists $R_V<\infty$ 
so that 
$$\limsup_{N\ra\infty}\mu^N_V(\mun_\bA( A^{2n}_i))\le 
(R_V)^n$$
for all $n\in\N$ and $i\in\{1,\cdots,m\}$. Here, $R_V$
is uniformly bounded by some $R_M$
when 
the quantities  $\left(a, V(0, 0,..,0), (\Da_i V(0, 0,..,0))_{1\le i\le m}\right)$
are bounded by $M$.

\item
$\mu^N_V[\mun_\bA]$ is tight 
and its limit points satisfy {\bf SD[V]}.

\item Take  $V=V_{\oo{t}}=\sum_{i=1}^n t_i q_i$ 
and let $U_a$ be the set of $t_i$'s
for which $V_{\oo{t}}$
satisfies (\ref{conv}) for a given 
$a<1$. 
For $\e>0$ small enough, when $(t_i)_{1\le i\le n} 
\in U_a\cap B(0,\e)$, $\mu^N_V[\mun_\bA]$ converges to
the unique solution to {\bf SD[V]}.

\item Assume that $U_a$ contains 
$\cup_{1\le i\le n}\{ (0,..,0, t_i,..,t_n), 0\le t_i\le \d\}$
for $\d$ small enough. Then, for $\e>0$ small enough,
for $\oo{t}\in U\cap B(0,\e)$, 
$$F_{V_{\oo{t}}}^N=N^{-2}
\log( {Z^N_{V_{\oo{t}}}})$$
converges as $N$ goes to infinity towards
$$F_{V_{\oo{t}}}=\sum_{\oo{k}\in \N^n\backslash (0,..,0)}
\prod_{1\le i\le n} \frac{(-t_i)^{k_i}}{k_i!}
\Ma((q_1,k_1),\cdots,(q_n,k_n)).$$
\end{itemize}
\end{theo}
{\bf Remark :} Observe that our hypothesis
is verified for  all quadratic interaction models
such as the Ising model, the $q$-Potts model
... etc ... as soon as the self potential 
of each matrix is convex.

\bigskip
\goodbreak
\noindent
\begin{dem}
We can assume without loss of generality that 
$a=0$ since otherwise we
just make a shift on the covariance 
of the matrices under  $\mu_N$.
The idea is to  use Brascamp-Lieb inequality (c.f Harge \cite{Har}
for recent improvements)
which  shows that since
$$f(\bA)= e^{-N\tr V(A_1,\cdots,A_m)}$$
is log-concave,
for all
convex function
$g$ on $(\R)^{mN^2}$
\begin{equation}\label{blieb}
\mu^N_V( g(\bA-\bM)) =\int g(\bA  -\bM
 ) \frac{f(\bA)\prod d\mu_N(A_i)}
{\int f(\bA)\prod d\mu_N(A_i)}
\le \int g(\bA)  \prod d\mu_N(A_i)
\end{equation}
with $$\bM =\int \bA d\mu^N_V.$$
Here $\bA$ denotes the set of entries of the matrices $(A_1,\cdots,A_m)$.
Let us apply (\ref{blieb}) with 
$g(\bA)=\tr(A^{2p}_k)$ which is convex by Klein's lemma.
Hence,
$$\mu^N_V\left( \tr((A_k-\E[A_k])^{2p}) \right)
\le \mu_N( \tr(A^{2p}_k))$$
where $\E[A_k](ij)=\mu^N_V(A_k(ij))$  for $1\leq i,j \leq N$
Since the right hand side is bounded by
$4^p$ as $N$ goes to infinity, we conclude 
that 
\begin{equation}
\label{jn}
\limsup_{N\ra\infty}\mu_V^N[\frac{1}{N}\tr((A_k-\E[A_k])^{2p}
)]\le 4^p.
\end{equation}
We now control $\E[A_k]$ uniformly.
Since the law of $A_k$ is
invariant by the action of the
unitary group, we deduce
that for all unitary matrix $U$,
\begin{equation}
\label{jn1}\E[ A_k]= \E[ UA_k U^*]= U\E[ A_k]U^*\Rightarrow 
\E[ A_k] =\mu_V^N( \mun_\bA(X_k)) I.
\end{equation}
We now bound $\mu_V^N( \mun_\bA(X_k))$ independently 
of $N$. Since $V$ is convex, there are real numbers $(\gamma_i)_{1\le
i\le m}$ and $c>-\infty$, $\gamma_i=\Da_i V(0,\cdots,0)$ and
$c=V(0,\cdots,0)$ so that
for all $N\in\N$ and all matrices $(A_1,\cdots,A_m)\in\Ha_N^m$,
$$\tr(V(A_1,\cdots,A_m))\ge \tr(\sum_{i=1}^m \gamma_i A_i +c).$$

By Jensen's inequality, we know that $Z_N^V\ge e^{-dN^2}$
for some $d<\infty$ and so Chebychev's
inequality implies that for all $y>0$, 
all $\l>0$, 
\begin{eqnarray*}
\mu_V^N\left( |\mun_\bA(X_k)|\ge y\right)
&\le& e^{(d-c)N^2-\l y N^2}[ \int e^{ -N\sum_{i=1}^m \gamma_i \tr(
A_i)
 + N \l \tr( A_k)
}\prod_{i=1}^m d\mu_N(A_i)\\
&&
+\int e^{ -N\sum_{i=1}^m \gamma_i \tr(
A_i)
 - N \l \tr( A_k)
}\prod_{i=1}^m d\mu_N(A_i)]\\
&\leq& 2e^{(d-c)N^2-\l y N^2} e^{\frac{N^2}{2} \sum_{i\neq k} 
\gamma_i^2 +\frac{N^2}{2}(\gamma_k+\l)^2}
\\
\end{eqnarray*}
Optimizing with respect to $\l$ shows that there exists $A<\infty$ 
so that
$$\mu_V^N\left( |\mun_\bA(X_k)|\ge y\right)\le e^{ AN^2 -\frac{N^2}{4} 
y^2}$$
and so $$\mu_V^N(|\mun_\bA(X_k)|)=\int 
\mu_V^N\left( |\mun_\bA(X_k)|\ge y\right) dy\le 4\sqrt{A}
+\int_{y\ge 4\sqrt{A}} e^{ -\frac{N^2}{4} 
(y^2-4A)} dy\le 8\sqrt{A}$$
where we assumed $N$ large enough in the last line.
Hence,
we have proved that
\begin{equation}
\label{jn2}\limsup_{N} |\mu_V^N(\mun_\bA(X_k))|< 8\sqrt{A}.
\end{equation}
Plugging this result in (\ref{jn}) and (\ref{jn1}) 
we obtain for all $p\geq 1$:
\begin{eqnarray*}
\limsup_{N\ra\infty}\mu_V^N[\mun_{\bA}((A_k)^{2p}
)]&\le&
2^{2p-1}\limsup_{N\ra\infty}\mu_V^N[\frac{1}{N}\tr ((A_k-\mu_V^N[A_k])^{2p}
)]\\
& &+2^{2p-1}\limsup_{N\ra\infty}(\mu_V^N( \frac{1}{N}\tr [A_k])^{2p}
)\\
&\le& 2^{2p-1} 4^p +2^{2p-1} (8\sqrt{A})^{2p}\le R_V^{2p}\\
\end{eqnarray*}
with $R_V=4(1+8\sqrt{A}).$
To prove the convergence of  $\mu_N^V[\mun_{\bA}]$, remember 
that $\mu_N^V[\mun_{\bA}]$ is tight for the $\Ca_{st}^m(\C)$-topology. To 
study its limit point, 
 recall $\int xe^{-x^2/2} f(x)dx=\int
f'(x)e^{-x^2/2} dx$ so that, for $P\in\Ca_{st}^m(\C)$,
\begin{eqnarray*}
\int \frac{1}{N}\tr(A_kP) d\mu_N^V(\bA)
&=&\frac{1}{2N^2}\sum_{ij}\int \partial_{A_k({ij})} (Pe^{-N \tr(V)})_{ji}
\prod d\mu_N(A_i)\\
&=&\frac{1}{2N^2}\sum_{ij}\int\left(\sum_{P=QX_kR}2Q_{ii}R_{jj}\right.\\
& &\left.-N\sum_{l=1}^n\sum_{q_l=QX_kR}t_l\sum_{h=1}^N2P_{ji}Q_{hj}R_{ih}\right)d\mu^N_V(\bA)\\
&=&\int
\left(\frac{1}{N^2}(\tr\otimes\tr)(D_k P)-
\frac{1}{N}\tr(\Da_k V P)\right)d\mu^N_V(\bA)
\end{eqnarray*}
which yields
$$\int \left( \mun_{\bA}((X_k+\Da_k V)P)- \mun_{\bA}\otimes  \mun_{\bA}
(D_kP)\right) d\mu_N^V(\bA)=0$$
Now, by convexity of $V$ we have concentration of $\mun_{\bA}$ under
$\mu_N^V$ (since log-Sobolev 
inequality is satisfied uniformly, according to 
Bakry-Emery criterion,
and that Herbst's 
argument therefore applies, see \cite{ane}, sections 6 and 7):
for all Lipschitz function $f$ on the entries
$$\mu_V^N\left(\bA: |f(\bA)-\mu_V^N(f)|\ge \d\right)\le e^{-\frac{\d^2}{2 ||f||_\La^2}}$$
where $||f||_\La$ is the Lipschitz constant of $f$.
Since for $P\in\Ca_{st}^m(\C)$, $\bA\ra \mun_\bA(P)$
is Lipschitz with constant of order $N^{-1}$ (see \cite{GZ1}),
we conclude that since $D_iP\in \Ca_{st}^m(\C)\otimes \Ca_{st}^m(\C)$,
for all $P\in\Ca_{st}^m(\C)$,
$$\lim_{N\ra\infty}\left|\int  \mun_{\bA}\otimes  \mun_{\bA}
(D_kP) d\mu_N^V(\bA)- \mu_N^V[\mun_{\bA}]\otimes \mu_N^V[\mun_{\bA}](D_kP)
 \right|=0.$$
Thus
$$\limsup_{N\ra\infty}
\left( \mu_V^N( \mun_{\bA}((X_k+\Da_k V)P))
- \mu_V^N[\mun_{\bA}]\otimes \mu_V^N[ \mun_{\bA}]
(D_iP)\right)=0$$
If $\tau$ is a limit point of $\mu_V^N[ \mun_{\bA}]$
for the weak $\Ca_{st}^m(\C)$-topology, we
can use the previous moment estimates to show that
even though $X_k+\Da_k V$ is a polynomial function,
$\mu_V^N( \mun_{\bA}((X_k+\Da_k V)P))$ converges 
along subsequences towards $\tau((X_k+\Da_k V)P))$,
and of course $\mu_V^N[\mun_{\bA}]\otimes \mu_V^N[ \mun_{\bA}]
(D_kP)$ converges towards $\tau\otimes \tau(D_kP)$.
Hence, we get that the limit points
of $\mu_V^N[\mun_{\bA}]$
satisfy the {\bf WSD[V]}.
By the previous moment estimate,
this limit points are compactly supported,
hence they satisfy  {\bf SD[V]}.

When $V=V_{\oo{t}}$, observe that $R_{{\oo{t}}}$
is uniformly bounded when $|t|\le M$
since $V_{\oo{t}}(0,\cdots,0)$ and $(\Da_i V_{\oo{t}}(0,\cdots,0))_{1\le
i\le m}$ depends continuously on $\oo{t}$.
Thus, the first point of the Theorem 
shows that the limit points
of  $\mu_{V_{\oo{t}}}^N[ \mun_{\bA}]$ 
are uniformly compactly supported. Hence, since also we have seen
 that they 
satisfy {\bf SD[$V_{\oo{t}}$]}, for $\oo{t}$ small enough, 
$\mun_\bA$ converges 
in expectation (and therefore almost surely
by concentration), to the unique solution to {\bf SD[$V_{\oo{t}}$]}.
The last point is now  a direct consequence of Theorem \ref{convfe}.

\end{dem}

Hence, we see here that 
convex potentials have uniformly compactly supported 
limit distributions so that we can apply 
the whole machinery. We strongly believe 
that this property extends to much more general
potentials. 
However, we shall
see in the next section that 
we can localize the
integral to make sure that
all limit points are uniformly compactly 
supported and still keep
the enumerative property, hence bypassing 
the issue of compactness.

\subsection{The uses of diverging integrals }\label{div}
In the domain 
of matrix models,
diverging integrals are often
considered.
For instance, if one
wants to consider
triangulations, one would like
to study the integral
$$Z_N(tx^3)=\int e^{tN\tr( M^3)} d\mu_N(M)$$
which is clearly infinite if $t$ is real.
The same kind of problem 
arises  in many 
other models (c.f.  the dually weighted graph model
\cite{KSW}).
However, we shall see below 
that at least as far as 
planar maps are concerned,
we can localize the integrals
to make sense of it, while keeping
its enumerative property.
Namely, let $V_{\oo{t}}=\sum t_i q_i$
and  let us consider the localized matrix
integrals given, for $L<\infty$,  by
$$Z_{V_{\oo{t}}}^{N,L}=\int_{||\bA||_\infty\le L}
e^{-N\tr(V_{\oo{t}}(\bA))} \prod d\mu_N(A_i)$$
and the associated Gibbs measure
$$\mu_{V_{\oo{t}}}^{N,L}(d\bA)=(Z_{V_{\oo{t}}}^{N,L})^{-1}
1_{||\bA||_\infty\le L}
e^{-N\tr(V_{\oo{t}}(\bA))} \prod d\mu_N(A_i).$$
Here, $||\bA||_\infty=\max_{1\le i\le m}||A_i||_\infty$
and $||A_i||_\infty$ denotes the spectral radius 
of the matrix $A_i$.

We  shall prove that

\begin{theo}\label{localization}
There exists  
$\e_0>0$ so that for $\e<\e_0$,
there exists $L_0(\e)$ and $L(\e)$, $L(\e)$ going to infinity  
and $L_0(\e)$ going to $2$ as $\e$
 goes to zero, so that for $\oo{t}\in B(0,\e)$,
for all $L\in [L_0(\e),L(\e)]$,
$$\lim_{N\ra\infty} \frac{1}{N^2}\log
Z_{V_{\oo{t}}}^{N,L}
=\sum_{\oo{k}\in \N^n\backslash (0,..,0)}
\prod_{1\le i\le n} \frac{(-t_i)^{k_i}}{k_i!}
\MM((q_1,k_1),\cdots,(q_n,k_n))$$

Moreover, under $\mu_{V_{\oo{t}}}^{N,L}$,
$\mun_\bA$ converges almost surely 
towards $\tau_{\oo{t}}$ described in Theorem \ref{existanaltheo}.
\end{theo}
This shows that, up to localization,
the first order asymptotics of matrix models 
give the right enumeration for any polynomials.
The diverging integrals often considered
in physics should be therefore thought to be conveniently
localized to keep their combinatorial
virtue, and are then as good as others.
In view of Lemma \ref{local}, this localization 
procedure should not 
damage the  rest of  the large $N$ expansion neither.
Note that when $m=1$, the localization amounts to restrict the
integral in the domain of strict convexity of
the potential, henceforth 
again avoiding all issues of
escaping eigenvalues. 

\nn
\begin{dem}
The proof is very close to that of Theorem \ref{mmSD}
except that we have to be careful
to make perturbations which do 
not change the constraint $||\bA||_\infty\le L$.
Let $i\in\{1,\cdots,m\}$ and consider the perturbation
$A_i\ra A_i +N^{-1} h(A_i)$ and $A_j\ra A_j$ for $j\neq i$
with a compactly supported function $h$ which vanishes 
on $[-R,R]^c$.
Then for $L>R$, for sufficiently large $N$, and $||A_i||_\infty\le L$,
$||A_i +N^{-1} h(A_i)||_\infty\le L$ so that
we see that the limit points of $\mun_\bA$
under the localized Gibbs measure
$$d\mu_V^{N,L}(A_1,\cdots,A_m)=(Z^L_N)^{-1} 1_{ ||A_i||_\infty\le L}
e^{-N\tr(V(\bA))}d\mu_N(A_1)\cdots d\mu_N(A_m)$$
satisfy for $i\in\{1,\cdots,m\}$,
\begin{equation}\label{eqcut}
\mu\otimes\mu(D_i h(X_i)) =\mu((\Da_i V+X_i) h(X_i)).
\end{equation}
These limit points 
are also laws of operators bounded
by $L$, but we shall see that in fact
this bound can be improved to become
independent of $L$ for good $L$'s.
We fix a limit point $\mu$ below; $\mu$ is a tracial state.
We proceed as before by taking $h(x)=P(\phi_R(x))$
with a polynomial test function $P$ and $\phi_R$ 
a smooth cutoff function with uniformly bounded 
derivative (say by $2$). Then,
$$\mu\otimes\mu(D_i h(X_i))=
\mu\otimes\mu(D_i P(\phi_R(X_i))\ts\phi'_R(X_i)\otimes 1)$$
so that if $P(x)=x^{2k+1}$, 
$$|\mu\otimes\mu(D_i h(X_i))|\le 2\sum_{p=0}^{2k} \mu(|\phi_R(X_i)|^p)
\mu(| \phi_R(X_i)|^{2k-p}).$$
Hence, we
get, for $i\in\{1,\cdots,m\}$,
$$\mu(X_i\phi_R(X_i)^{2k+1})\le 2\sum_{p=0}^{2k} \mu(|\phi_R(X_i)|^p)
\mu( |\phi_R(X_i)|^{2k-p})
+|t|\sum_{j=0}^{n}\mu(|\Da_i q_j(\bX)||\phi_R(X_i)|^{2k+1})$$
Note that we can bound above 
the last term by H\"older's inequality
so that
$$\sum_{j=0}^{n}\mu(|\Da_i q_j(\bX)||\phi_R(X_i)|^{2k+1})
\le C\max_{1\le j\le m}\{\mu( |X_j|^{2k+D-1}),
\mu( |X_j|^{2k})\}$$
where we assumed $|\phi_R(x)|\le |x|$.
Taking $\phi_R(x)= x$ when $|x|\le R/2$ and 
$\phi_R(x)=0$ when $|x|\ge R$, $\phi_R$ linear 
in between,
we can now use 
monotone convergence theorem (letting $R$ 
going to infinity) to obtain
$$\max_{1\le j\le m}\mu(X_j^{2(k+1)})\le
2\sum_{p=0}^{2k} \max_{1\le j\le m}\mu(|X_j|^p)
\max_{1\le j\le m}\mu( |X_j|^{2k-p})
+C|t|\max_{1\le j\le m}(\mu(X_j^{2k})+\mu(|X_j|^{D+2k-1} ))$$
Noting that $\mu(|X_j|^{D+2k-1} )\le L^{D-1}\mu(|X_j|^{2k})$
since under $\mu$ the operators are uniformly bounded by $L$,
we can improve this uniform bound as follows.
We make the induction hypothesis that
$$B_k:=\max_{1\le j\le m}\mu(|X_j|^k)\le C_k A^k$$
for $k\le 2r$,
with $C_k$ the Catalan numbers and some $R>1$.
Then, by H\"older's inequality,
$$ (B_{2r+1})^{\frac{2r+2}{2r+1}}\le B_{2r+2}\le  2 C_{2r+1}A^{2r} +C|t|(1+L^{D-1}) A^{2r}C_{2r}.$$
We can assume without loss of generality that $B_{2r+1}\ge 1$ so that
we get 
$$B_{2r+2}\vee B_{2r+1}\le C_{2r+1} A^{2r+1}\left( 2A^{-1}  +C|t|L^D\right)$$
and so we can continue our induction if $L$ is not too large ($L$ at most of order
$ (2 C|t|)^{-\frac 1 D}$)
so that there is $A>0$ so that $2A^{-1}+C|t|L^D\le 1$. 
Hence  all limit points are
non-commutative laws
of operators with norms uniformly bounded by $4A\ll L$.
Thus if $R>4a$ we conclude by (\ref{eqcut}) that all limit points satisfy {\bf SD[$V_{\oo{t}}$]}.
Moreover, Theorem \ref{unique} and \ref{existanaltheo} apply
to show that there is a unique solution to this
equation which are laws of operators bounded by $R$ (provided $R$,
and so $L$, is neither too big ($L\le O(|t|^{-\D^{-1}})$,
for uniqueness), 
nor too small ($L\ge L_0( |t|)\ge L_0(0)=2$ 
for existence. Note here that $L_0( |t|)$
is the smallest real number such that $|\tau_{\oo{t}}(P)|\le L_0(
|t|)^{deg(P)}$ for all monomial $P$, which converges to $2$ as 
$\oo{t}$ goes to zero by Theorem \ref{existanaltheo}.
Hence, $\mun_\bA$ has a unique limit point, $\tau_{\oo{t}}$,
and thus converges towards it.

The  formula of the free energy is then 
derived as in Theorem \ref{convfe}
since $L$ is fixed independently of $\oo{t}$ 
small enough.

\end{dem}
Let us remark
that if we define, following Voiculescu \cite{Vo1},
a microstates $\G(\mu, n,N, \e)$
for $\mu\in \Ma_{1}^{(m)}
$,
$n\in\N$,
$N\in\N$, $\e>0$, as the set of  matrices $A_1,..,A_m$ of $\Ha_N^m$ such that
\begin{equation}\label{microstates}
|\mu(\bX_{i_1}..\bX_{i_p})
-\trn(\bA_{i_1}..\bA_{i_p})|<\e
\end{equation}
for any  $1\le p\le n$,
$i_1,..,i_p\in\{1,..,m\}^p$, then we have 

\begin{lem}\label{local}
For $L$ big enough, $V=V_{\oo{t}}$ 
and  $\oo{t}$ small enough,
$$\lim_{N\ra\infty}
\frac{1}{ N^2}\log \int_{||\bA||_\infty\le L}
e^{-N\tr(V(\bA))}d\mu_N(A_1)\cdots d\mu_N(A_m)
$$
$$
=\lim_{\e\ra 0, n\ra\infty}\lim_{N\ra\infty} 
\frac{1}{ N^2}\log  \int_{\G(\tau_V, n,N, \e)\cap ||\bA||_\infty\le L}
e^{-N\tr(V(\bA))}d\mu_N(A_1)\cdots d\mu_N(A_m)$$
$$
=\lim_{\e\ra 0, n\ra\infty}\lim_{N\ra\infty}
 \frac{1}{ N^2}\log  \int_{\G(\tau_V, n,N, \e)}
e^{-N\tr(V(\bA))}d\mu_N(A_1)\cdots d\mu_N(A_m)$$

\end{lem}
\begin{dem}
The first equality is a direct consequence 
of the previous Theorem since it is equivalent
to the fact that $\mu_{V_{\oo{t}}}^{N,L}( \G(\tau_V, n,N, \e))$
goes to one. The second 
comes from the fact that for $n$ greater than
the degree of 
$V$, 
$$\lim_{\e\ra 0, n\ra\infty}\lim_{N\ra\infty} 
\frac{1}{ N^2}\log  \int_{\G(\tau_V, n,N, \e)\cap ||\bA||_\infty\le L}
e^{-N\tr(V(\bA))}d\mu_N(A_1)\cdots d\mu_N(A_m)$$
$$= -\tau_V( V) 
+ \lim_{\e\ra 0, n\ra\infty}\lim_{N\ra\infty} 
\frac{1}{ N^2}\log  \mu_N^{\ot m}\left( 
\G(\tau_V, n,N, \e)\cap ||\bA||_\infty\le L\right)$$
$$= -\tau_V( V) 
+ \lim_{\e\ra 0, n\ra\infty}\lim_{N\ra\infty} 
\frac{1}{ N^2}\log  \mu_N^{\ot m}\left( 
\G(\tau_V, n,N, \e)\right)$$
where we used in the last equality the result of \cite{bebe},
which hold when $\tau_V$ is the law
of bounded operators with norms strictly smaller than $L$.
\end{dem}
Therefore, the localization should not affect 
the full expansion of the integral
since second order asymptotics are usually obtained 
first by a localization on a microstates in order to use 
precise Laplace method's.

As a corollary, we also deduce that
for all $V_{\oo{t}}$ with ${\oo{t}}$ 
small enough,  the limits of empirical
distributions of matrices
given by localized matrix models 
provide solutions of {\bf SD[$V_{\oo{t}}$]}.
Since these limits have to be 
tracial states, we deduce that

\begin{cor}  The solutions compactly supported of {\bf SD[$V_{\oo{t}}$]} 
are tracial states when $\oo{t}$ is sufficiently small.

\end{cor}
Note that if $(P_i)_{1\le i\le m}$ 
is the conjugate variable
of a tracial state, Voiculescu \cite{VoCyc}
have shown that $P_i=\Da_i P$
for $1\le i\le m$ and some polynomial $P$.
This fact should be compared with our graphical
interpretation which works only because $P_i$
is a cyclic derivative.
\section{Applications  to free entropy}
Let us recall that Voiculescu's microstates
entropy is defined, for $\tau\in\cup_R \Ma^m_R$,
by
$$\chi(\tau)
=\lim_{\e\ra 0, n\ra\infty\atop
L\ra\infty}\limsup_{N\ra\infty}\frac{1}{N^2}
\log \mu_N^{\ot m}\left(\Gamma(\tau,n,N,\e)\cap||\bA||_\infty\le
L\right)$$
with $\Gamma(\tau,n,\e,N)$ the microstates defined in 
(\ref{microstates}).
Note that the original definition 
of Voiculescu is not with respect to the Gaussian measure,
but with respect to the Lebesgue measure.
However, both definitions only differ by a 
quadratic term (see \cite{CDG2}).
It is an (important) open problem 
whether in general one can replace 
the limsup by a liminf in the definition of
$\chi$. However, from the previous considerations,
we can see the following

\begin{theo}
\label{freeentr}
Let $n\in \N$ and $(q_i)_{1\le i\le n}$
be  monomials in $m$ non-commutative
variables $\bX=(X_1,\cdots,X_m)$. Let $V_{\oo{t}}(\bX)=
\sum_{i=1}^n t_i (q_i(\bX)+q_i^*(\bX))$. By Theorem \ref{existanaltheo},
we know that there exists $\e>0$ so that 
for $|t|<\e$,
there exists a unique solution $\tau_{\oo{t}}$
to {\bf SD[$V_{\oo{t}}$]}. Then,
also for $|t|\le \e$,
$$\chi(\tau_{\oo{t}})
=\lim_{\e\ra 0, n\ra\infty\atop
L\ra\infty}\liminf_{N\ra\infty}\frac{1}{N^2}
\log \mu_N^{\ot m}\left(\Gamma(\tau,n,N,\e)\cap \{||\bA||_\infty\le L\}
\right).$$
Moreover, 
$$\chi(\tau_{\oo{t}})=-\sum_{k\in\N^n\backslash (0,\cdots,0)}
\prod_{i=1}^n \frac{(-t_i)^{k_i}}{k_i!} \left(\sum_{j=1}^n k_j-1\right) \Ma((q_1,k_1),\cdots,(q_n, k_n)).$$

\end{theo}

\nn
{\bf remark:} In particular, we see 
as could be expected that $\chi(\tau_{\oo{t}})>-\infty$
and so all the  solutions to Schwinger-Dyson's equations
for small parameters are laws of von Neumann
algebras which are isomorphic
to the free group factor.

\begin{dem}
In fact,
\begin{eqnarray*}
\chi(\tau_{\oo{t}})&=&\lim_{\e\ra 0, n\ra\infty\atop
L\ra\infty}\limsup_{N\ra\infty}\frac{1}{N^2}
\log 
\int_{ \Gamma(\tau_{\oo{t}},n,N,\e)\cap \{||\bA||_\infty\le L\}} e^{ N\tr(V_{\oo{t}}(\bA))
-N\tr(V_{\oo{t}}(\bA))} d \mu_N^{\ot m}(\bA)\\
&=&\tau_{\oo{t}}(V_{\oo{t}})+
\lim_{\e\ra 0, n\ra\infty\atop
L\ra\infty}\limsup_{N\ra\infty}\frac{1}{N^2}
\log 
\int_{ \Gamma(\tau_{\oo{t}},n,N,\e)\cap \{||\bA||_\infty\le L\}} e^{
-N\tr(V_{\oo{t}}(\bA))} d \mu_N^{\ot m}(\bA)\\
&\le& \tau_{\oo{t}}(V_{\oo{t}})
+ F_{\oo{t}}\\
\end{eqnarray*}
with
$$F_{\oo{t}}=\lim_{ 
L\ra\infty}\limsup_{N\ra\infty}\frac{1}{N^2}
\log 
\int_{||\bA||_\infty\le L}  e^{
-N\tr(V_{\oo{t}}(\bA))} d \mu_N^{\ot m}(\bA).$$
On the other hand,
\begin{eqnarray*}
&&\lim_{\e\ra 0, n\ra\infty\atop
L\ra\infty}\liminf_{N\ra\infty}\frac{1}{N^2}
\log \mu_N^{\ot m}\left(\Gamma(\tau_{\oo{t}},n,N,\e)\cap||\bA||_\infty\le L
\right)\\
&=&\tau_{\oo{t}}(V_{\oo{t}})+
\lim_{\e\ra 0, n\ra\infty\atop
L\ra\infty}\liminf_{N\ra\infty}\frac{1}{N^2}
\log 
\int_{ \Gamma(\tau_{\oo{t}},n,N,\e)\cap \{||\bA||_\infty\le L} e^{ 
-N\tr(V_{\oo{t}}(\bA))} d \mu_N^{\ot m}(\bA)\\
&=&\tau_{\oo{t}}(V_{\oo{t}})+ F_{\oo{t}}
+\lim_{\e\ra 0, n\ra\infty\atop
L\ra\infty}\liminf_{N\ra\infty}\frac{1}{N^2}
\log \mu_{V_{\oo{t}}}^{N,L}\left(  \Gamma(\tau_{\oo{t}},n,N,\e)
\right)\\
&=& \tau_{\oo{t}}(V_{\oo{t}})+ F_{\oo{t}}\\
\end{eqnarray*}
where we used in the last term Theorem \ref{localization}
which implies
$$\lim_{N\ra\infty} \mu_{V_{\oo{t}}}^{N,L}\left(  \Gamma(\tau_{\oo{t}},n,N,\e)
\right)=1$$
for all $\e>0,n\in\N$ and $L$ large enough.
Thus, we see that $\chi$ is equal to its liminf definition
and moreover
$$ \chi(\tau_{\oo{t}})=\tau_{\oo{t}}(V_{\oo{t}})+ F_{\oo{t}}.$$
Now, by Theorems \ref{localization}
and \ref{existanaltheo}, 
$$F_{\oo{t}}=\sum_{\oo{k}\in \N^n\backslash (0,..,0)}
\prod_{1\le i\le n} \frac{(-t_i)^{k_i}}{k_i!}
\Ma((q_1,k_1),\cdots,(q_n,k_n))$$
whereas
$$\tau_{\oo{t}}(V_{\oo{t}})=\sum_{i=1}^n t_i 
\sum_{  k_j\in \N,\atop 1\leq j\leq n} \prod_{1\le j\le n} 
\frac{(-t_j)^{k_j}}{k_j!}\Ma((q_1,k_1),\cdots, (q_{i-1},k_{i-1}),
(q_i,k_i+1),(q_{i+1},k_{i+1}),\cdots, (q_n,k_n))$$
from which the formula for $\chi(\tau_{\oo{t}})$
is easily derived.
\end{dem}

\section{Applications to the combinatorics of planar maps}

For the sake of completeness,  we summarize in this last section,
the results of a few 
papers devoted to the enumeration
of planar maps, either by a 
combinatorial approach
or by a matrix model
approach.

\subsection{The one matrix case}
We now consider the simpler case $m=1$ where we only have one matrix.
Let $V_{\oo{t}}(A)=\sum_{i=1}^{2D}t_iA^i$ with $t_{2D}>0$ a polynomial
potential with an even 
leading power.
Then it has been proven in \cite{BAG} 
Theorem 5.2 that the empirical measure satisfies a large deviation principle:

\begin{theo}
Let $$J(\mu)=\int \left(\frac{x^2}{2}+
V_{\oo{t}}(x)\right)d\mu(x)-\int\int\log|x-y|d\mu(y)
d\mu(x)$$
and $$I(\mu)=J(\mu)-\inf_{\nu\in P(\RR)}J(\nu)$$
then the sequence of empirical measure $\mun$ satisfies a large 
deviation principle
in the scale $N^2$ with good rate function $I$. 
Moreover,
the minimum of $I$ is reached at a unique probability measure 
$\mu_{\oo{t}}$ so that
$$\frac{x^2}{2}+V_{\oo{t}}(x)-
2\int\log|y-x|d\mu_{\oo{t}}(y)=
C_{\oo{t}},\,\,\,\mu_{\oo{t}} a.s.$$
with a finite constant $C_{\oo{t}}$,
and where the left hand side dominates 
the right hand side on the whole
real line. 
\end{theo}

One can notice that differentiating in $x$ the last equation,
we recover the Schwinger Dyson's equation. It is not sufficient
in general
to determine the solution uniquely; one  
need the  inequality on the whole real line
to fix the support of the solution.

These questions have also been investigated with 
the method of orthogonal polynomials which give
 a rather sharp description 
of the limit measure and emphasizes
 a structure similar to the semi-circular law.
More precisely  Theorem 3.1 in \cite{EML} gives:

\begin{theo}\label{deift}
There exists $t>0$ and $\gamma>0$ such that if for all $i$, $|t_i|<t$
 and $t_{2D}>\gamma\sum_{i<2D}t_i$ then $\mu_{\oo{t}}$ is absolutely
continuous with density $\Psi_{\oo{t}}$ of the form:
$$\Psi_{\oo{t}}(x)=\frac{1}{2\pi}1_{[a,b]}(x)\sqrt{(x-a)(x-b)}h(x)$$
with
$$h(z)=\int_{C(z,R)}\frac{V_{\oo{t}}'(s)}{\sqrt{(s-a)(s-b)}}\frac{ds}{s-z}$$
where $R$ is such that $a,b\in C(z,R)$.
Besides, the boundaries $a$ and $b$ can be find by the equations:
$$\int_a^b\frac{V_{\oo{t}}'(s)}{\sqrt{(s-a)(b-s)}}ds=0$$
$$\int_a^b\frac{sV_{\oo{t}}'(s)}{\sqrt{(s-a)(b-s)}}ds=2\pi$$
\end{theo}

We now look at combinatorics of the Schwinger-Dyson's 
equation with one variable, for $V_{\oo{t}}(x)=\sum_{i=1}^{2D}t_ix^i$.
Remember that from Theorem \ref{existanaltheo}, 
$\mu_{\oo{t}}$ can be seen
as the generating function of graphs 
counted by the numbers of stars of valence $i$:
$$\mu_{\oo{t}}(x^p)=\sum_{k_1,\cdots,k_{2D}\in\NN}
\prod_{i=1}^{2D}\frac{(-t_i)^{k_i}}{k_i!}\MM_0((P,1),
\{(X^i,k_i)\}_{1\leq i\leq 2D}).$$
Hence, Theorem \ref{deift} allows 
to estimate the numbers of
one colour planar maps.
A more direct combinatorial approach 
can be developed by considering for instance 
the dual of those graphs.
The dual of a graph is simply obtained by 
replacing each face by a star 
and each edge by a transverse edge which link the
 two stars which come from
the face adjacent to the edge. In that operation 
each star is replaced by a face of the same valence. 
As we work on the sphere we can decide that the face
 which comes from the star $X^p$ is the external face.

So $\mu_{\oo{t}}(X^{p+1})$ is also the generating 
function of connected planar graphs with an external face
of valence $p+1$ and enumerated by the number of faces of a given valence.
Those objects are classical ones in combinatorics and we can follow \cite{Tu}
to find an equation on these generating functions.
The idea is to try to cut the first edge of the external face, then two 
cases may occur: either the graph is disconnected and we obtain two graphs or 
it isn't disconnected and the external face has grown.
This two cases corresponds in the dual graph to the fact that the first branch 
of the root is a loop or not which is exactly what we use to build our 
combinatorial interpretation
so that we can retrieve the Schwinger Dyson's equation from this fact.
Just by using the equation given by this decomposition and 
some algebraic tools
combinatoricians have solved some  models.
For example \cite{BC} gives an equation 
on the generating function $M(u,v)$ of maps
 whose internal faces have degree living 
in a fixed set $\DD\subset \NN$ and enumerated by 
their number
of edges and the degree of the external face.
To translate this in our framework, one can consider 
for a finite $\DD$ with an even maximal element, 
$$V_{\oo{t}}(X)=\sum_{d\in \DD}t_dX^d$$
Then under this potential, for small $t$, the limit measure $\mu_{\oo{t}}$ 
will satisfy our combinatorial interpretation. Then 
$$M(u,v)=\sum_{p\in\NN}\mu_{(-u^{\frac{d}{2}})_{d\in \DD}}(X^p)v^p$$
Now Theorem 1 of \cite{BC} states:

\begin{theo}
For a serie $F(z)=\sum_ia_iz^i$ we will note $[z^i]F(z)$ 
the $i^{th}$ coefficient $a_i$. Then there exists 
a unique power serie $R$  satisfying
$$R=1-4R_1v-4R_2v^2$$
with 
$$R_1=\frac{u}{2}\sum_{i\in D}[v^{i-1}](R^{\frac{1}{2}})
\mbox{ and }R_2=\frac{u}{2}\sum_{i\in D}[v^{i}](R^{\frac{1}{2}})
+u-3R_1^2.$$
The number $m_n$ of maps with $n$ edges such that every degree of internal face lies in $D$ is then
$$m_n=[u^n]\frac{(R_2(u)+R_1(u)^2)(R_2(u)+9R_1(u)^2)}{(n+1)u^2}$$
\end{theo}
The techniques to prove these results are most 
often purely algebraic. The main difference 
in nature than we could meet between
the approaches by matrix models
or by combinatorics to these enumerations
is that
the first provides for free additional structure;
it shows 
that these enumerations  can be expressed 
in terms of a probability measure. This
point generalizes to any number of colors
where the enumeration can be expressed in terms
of tracial states. 
One may hope that this positivity condition 
could help in solving this combinatorics problems.

\subsection{Ising model on random graphs}\label{Isingsec}
This model is defined 
by $m=2$ and 
$$V(A,B)=V_{Ising}(A,B)=-cAB +V_1(A)+V_2(B).$$
In the sequel, we denote in short $A$ for $X_1$
and $B$ for $X_2$.
It is clear that for $|c|<1$, $V$ is a convex potential
as defined in \eqref{conv}  if $V_1,V_2$ 
are convex (write $-2AB=(A-B)^2 -A^2-B^2$
or $2AB=(A+B)^2 -A^2-B^2$ to see that up to a quadratic term
$2^{-1}|c| A^2+2^{-1} |c| B^2$, $V$ is convex)
Hence we deduce from Theorem \ref{convex} that

\begin{cor}\label{Ising}
For $c\in \R$ and  $V_i(x)=\sum_{j=1}^{D} t^i_j x^{2j}$, $i=1,2$, 
set
 $V_{\oo{t},c}(A,B)=-cAB +V_1(A)+V_2(B).$
Let,  for $\d>0$,  $U_\d=\cap_{i,j}\{ 0\le t_j^i\le \d\}\cap\{|c|<1-\d\}$.
Then, for $\d>0$ small enough and $(\oo{t},c)\in U_\d$,
$\mu_V^N( \mun_{\bA})$ converges towards
the solution $\mu_{\oo{t},c}$ of {\bf SD[$V_{\oo{t},c}$]} 
as $N$ goes to infinity. Moreover 

$$\mu_{\oo{t},c}(P)=\sum_{{\oo{k}}\in{\NN}^{2D}\atop r\in\NN}\prod_{i,j}\frac{(-t^i_j)^{k^i_j}}{k^i_j!}\frac{c^r}{r!}\MM_0((P,1),(A^{2j},k^1_j)_{1\leq j\leq D},(B^{2j},k^2_j)_{1\leq j\leq D},(AB,r)),$$

and 

$$F(\oo{t},c)-F(\oo{t},0)=
\sum_{{\oo{k}}\in{\NN}^{2D}
\atop r\geq 1}\prod_{i,j}\frac{(-t^i_j)^{k^i_j}}{k^i_j!}
\frac{c^r}{r!}\MM_0((A^{2j},k^1_j)_{1\leq j\leq D},
(B^{2j},k^2_j)_{1\leq j\leq D},(AB,r)).$$

\end{cor}
\noindent
{\bf Remark:}
Note that we took potentials $V_1$ and $V_2$ 
as polynomials with even powers to guarantee
our convexity relation but this 
condition could easily been relaxed 
by taking more sophisticated domains
than $U_\delta$ in which the polynomials would remain convex.

\noindent
\begin{dem} 
This result is a consequence
of Theorem \ref{main}, \ref{convex} and \ref{convfe}. Note here that
the control on $\mu^N_V(N^{-1}\tr(AB))$ assumed in Theorem \ref{convfe} 
is satisfied due to Theorem \ref{convex}
which provides a uniform bound when $|c|<\xi$
for $\xi<1$. 
\end{dem}

According to the graphical interpretation, the limiting measure 
is linked to planar maps with stars 
whose type are the monomial of $V_1, V_2$ 
and stars of type $AB$.
Those maps are very close from Ising configuration on planar graphs except that two stars of type $AB$ can be linked together.
For integers $(k^i_j)_{i\in\{A,B\},1\leq i\leq D}$, define
\begin{eqnarray*}
\II(\{k^i_j\},r,P)&=\sharp\{&\textrm{planar maps with $k^i_j$ stars of color $i$ and degree $2j$,}\\
& &\textrm{ one star of type $P$ (if $P\neq 0$) and $r$ stars of type $AB$ }\\
& &\textrm{ such that there's no link between any of the $r$ $AB$-stars. }\}
\end{eqnarray*}
and its rooted counterpart:

\begin{eqnarray*}
\JJ(\{k^i_j\},r,P)&=\sharp\{&\textrm{rooted planar maps with $k^i_j$ stars of color $i$ and degree $2j$,}\\
& &\textrm{ one star of type $P$ wich is the root and $r$ stars of type $AB$ }\\
& &\textrm{ such that there's no link between any of the $r$ $AB$-stars. }\}
\end{eqnarray*}

There's a relation between these quantities similar to (\ref{mapformula}):
\begin{equation}\label{isingmapformula}\II(\{k^i_j\},r,P)=\JJ(\{k^i_j\},r,P)r!\Pi_{i,j}k^i_j!(2j)^{k^i_j}\end{equation}

We can now relate these numbers to our limit measure:

\begin{prop}\label{Ising2}
Let $\mu_{\oo{t},c}$ be as in \ref{Ising}, then on its radius of convergence,

$$\mu_{\oo{t},c}(P)=\left(\frac{1}{1-c^2}\right)^{\frac{\textrm{deg }P}{2}}\sum_{k^i_j\in \NN^{2D}\atop r\in \NN}\prod_{i,j}\frac{1}{k^i_j!}\left(\frac{-t^i_j}{(1-c^2)^j}\right)^{k^i_j}\frac{c^r}{r!}\II(\{k^i_j\},r,P)$$
and
$$F(\oo{t},c)-F(\oo{t},0)=\frac{1}{1-c^2}\sum_{k^i_j\in \NN, 
i\in\{1,2\}, j\in \{1,D\}, r\geq 1}\prod_{i,j}\frac{1}{k^i_j!}\left(\frac{-t^i_j}{(1-c^2)^j}\right)^{k^i_j}\frac{c^r}{r!}\II(\{k^i_j\},r,0)$$
\end{prop}
\begin{dem} 
First we define a projection $\pi$ from rooted maps to rooted Ising graph such that if $M$ is a map $\pi(M)$ is obtained by deleting pairs of $AB$ stars which are glued.
We now apply Corollary \ref{Ising}, and translate its result in term of rooted diagrams using (\ref{mapformula}):

$$\mu_{\oo{t},c}(P)=\sum_{{\oo{k}}\in{\NN}^{2D}\atop r\in\NN}\prod_{i,j}(-2jt^i_j)^{k^i_j}c^r\DD_0(P,(A^{2j},k^1_j)_{1\leq j\leq D},(B^{2j},k^2_j)_{1\leq j\leq D},(AB,r))$$

All the maps $M$ appearing in that sum are such that $\pi(M)$
is an Ising graph rooted at a star of type $P$.
For a fixed Ising graph $G$ we must find the contribution in that sum 
of $\pi^{<-1>}(G)$.
But we can construct every graph in that set by adding pairs of
stars $AB$ on the edges of $G$. The numbers of edges 
of $G$ is $e_G=\frac{deg P}{2}+\sum_{i,j}2jk^i_j$ so that to get the whole 
contribution of $\pi^{<-1>}(G)$ we have to multiply the contribution of $G$ by
$$\sum_{a_1,\cdots,a_{e_G}\in\NN}c^{2\sum a_i}=
\left(\frac{1}{1-c^2}\right)^{\frac{deg P}{2}+\sum_{i,j}2jk^i_j}.$$
In that sum, $a_i$ stands for the number of pairs of $AB$ stars added on the $i^{th}$ edge.
Summing on every graphs, we obtain:

$$\mu_{\oo{t},c}(P)=\left(\frac{1}{1-c^2}\right)^{\frac{\textrm{deg }P}{2}}\sum_{{\oo{k}}\in{\NN}^{2D}\atop r\in\NN}\prod_{i,j}\left(\frac{-2jt^i_j}{(1-c^2)^{2j}}\right)^{k^i_j}c^r\JJ(\{k^i_j\},r,P)$$
and the result follows by using (\ref{isingmapformula}).
\par

The second point can be proven by proceeding in the same way.
\end{dem}

In the rest of this section,  we 
 compare  a few  different 
 approaches to solve
the enumeration problem  of 
the Ising model.
In short, let us
emphasize that, for the time 
being, combinatorial
and orthogonal polynomials approaches
give the more complete 
and explicit results.
However, these techniques are
still limited to very few models.
The Schwinger-Dyson's
equation or the large deviation
approaches can be developed
for a much wider range of models
(such as $q$-Potts, induced  QCD etc).
However, it seems to
us that these arguments
still need some mathematical
efforts to provide
as transparent and powerful
results (namely for the first
a mathematical  study of the so-called master-loop
 equations, and for the
second
a clear understanding 
of the relations between 
complex Burgers equations
and the master-loop equations). 
A striking difference between 
the combinatorial
 and the matrix
model approaches seems
to reside in the fact that
matrix models provide 
for free information
on the structure of 
the generating 
function of the number
of planar maps, for instance
as the Stieljes transform
of a probability measure 
with connected support.

\subsubsection{Orthogonal polynomial approach}
Here we take $V_1=V_2=(g/4) x^4$.
By using orthogonal 
polynomials techniques, it was proved by Mehta 
\cite{Me2}
that the corresponding free energy $F_{g,c}$ satisfies

$$F_{g,c}-F_{0,c}=\int_0^1 (1-x) [\log f(x)-\log \frac{cx}{2(1-c^2)}] dx$$
with $f(x)=f_{g,c}$ solution to the algebraic equation
$$f(x)\{ (1-6\frac{g}{c} f(x))^{-2} -c^2\}
+12 g^2 f^3(x) -\frac{1}{2} cx =0$$
and the root to be taken equals $2^{-1} cx(1-c^2)^{-1}$ 
when $g=0$.
\par
Starting from there, a simpler expression as been derived  in \cite{BK2} (equation (16), (17) with $h=z/g$):
\begin{eqnarray*}
F_{z,c}&=&\frac{1}{2}\ln h(z)+\frac{h^2(z)}{2}\left(\frac{z-1}{2(3z-1)^3}+c^2\frac{z+1}{3z-1}+\frac{c^4}{2}(3z^4-3z^2+1)\right)\\
& &-h(z)(\frac{1}{3z-1}+c^2(1-z^2))+\frac{1}{2}\ln(1-z^2)+\frac{3}{4}
\end{eqnarray*}
with
\begin{equation}\label{algebricequation1}
h(z)=\frac{(1-3z)^2}{1-c^2(1-3z)^2(1-3z^2)}
\end{equation}
Hence, by the preceding, Mehta's result
gives a formula for the generating 
function of $\JJ$ in the quadrangulation case.
However, it does not a priori gives
the limiting spectral measures of
the matrices. Moreover, this strategy could
be only developed completely and rigorously
for the Ising model
and the matrix coupled in chain model \cite{CMM}.

\par

\subsection{Direct combinatorial
approach}

We can also relate this result to the work of Bousquet-Melou 
and Schaeffer  \cite{BM-S}. Their approach is purely
combinatorial; they use bijection
with well labeled trees (whose 
generating functions are well understood)
to obtain algebraic equations
for the generating functions
of the Ising model.
Let $I(X,Y,u)$ be the generating function of the Ising model 
on quasi-tetravalent graphs, (i.e. tetravalent graphs except
 for the root which is bivalent and black) where $X$ (resp. $Y$) 
counts the black (resp. white) tetravalent stars and $u$ the 
bicolored edges:
\begin{eqnarray*}
I(X,Y,u)=& \sum_{m,n,r\in \NN}X^mY^nu^r\sharp\{&\textrm{quasi
-tetravalent maps with $m$ tetravalent black stars,}\\
& &\textrm{ $n$ tetravalent white stars and $r$ bi-colored edges}\}.
\end{eqnarray*}
If $P(x,y,u)$ is the solution to the algebraic equation:
\begin{equation}\label{algebricequation2}
P=1+3xyP^3+\frac{P(1+3xP)(1+3yP)}{u^2(1-9xyP^2)^2}
\end{equation}
Then, by \cite{BM-S}, Proposition 1 p.4, $I$ can be written in
 function of $P(x,y,u)$ with $x=X(u-\frac{1}{u})^2$ and 
$y=Y(u-\frac{1}{u})^2$ as
$$I(X,Y,u)=(1-u^{-2})\left(xP^3+\frac{P(1-3xP-2xP^2-6xyP^3)}{1-9xyP^2}
-\frac{yu^{-2}P^3(1+3xP)^3}{(1-9xyP^2)^3}\right).$$
On the other hand, according to  Proposition \ref{Ising2},
if $V=tA^4+uB^4-cAB$ and $\mu_{t,u,c}$ is the associated 
limit measure then on its domain of convergence,

 $$I(X,Y,u)=(1-u^2)\mu_{{X}{(1-u^2)^{-2}},{Y}{(1-u^2)^{-2}},u}(A^2).$$
If we make the following change of variable in 
(\ref{algebricequation2}):

$$x=y=\frac{-z}{3c^2h(z/3)}, P=-c^2h(z/3), u=c$$
then we find (\ref{algebricequation1}).
Hence, a combinatorial
approach can be developed to solve
the problem of the enumeration of 
planar maps of the Ising model,
a strategy which requires some 
 combinatorial
insight. The next approach we
present, developed in particular by
Staudacher, Kazakov and  Eynard,
is a direct analysis of the {\bf SD[V]}
equations. It is a purely analytical 
and rather robust strategy.

\subsection{Direct study of the {\bf SD[$V_{Ising}$]} equations}
Here, the analysis 
is based on Theorem \ref{convex}
which asserts that if $V_1,V_2$
are convex, for small parameters, $\mun_{A,B}$
converges almost surely
towards the solution $\mu_{\oo{t},c}$ of 
{\bf SD[$V_{\oo{t},c}$]}
which is a generating function 
for the enumeration of maps. 
Hereafter we take $c=1$ up to a rescaling $\bar x=\sqrt{c} x$,
$ \bar y=\sqrt{c} x$, $V_1(x)=\bar V_1(\bar x)$, 
$\mu_{\oo{t}}(P(A,X_2))=\mu_{\oo{t},1}
(P(\sqrt{c}^{-1} A,\sqrt{c}^{-1} X_2))$.
Following Eynard \cite{eynard},
we shall analyze the solutions 
of the Schwinger Dyson's equation.
Observe that the following considerations
hold for any range of parameters,
not only small parameters.
For large parameters, we do not know
that the Schwinger Dyson's equation
has a unique solution but we still
know that any limit point of the empirical 
measure of the random matrices
still satisfies it. In the next section, we shall
see that for the Ising model
and any range of parameters,
there is a unique such limit point,
and it will therefore enjoy the properties 
described below. We here summarize the
main result, as found in Eynard \cite{eynard}.
Let $\mu_{\oo{t}}$ be a solution
of  {\bf SD[$V_{Ising}$]}

$$\mu_{\oo{t}}((W_1'(A)-B)P)=\mu_{\oo{t}}\otimes\mu_{\oo{t}}(D_AP),$$
$$\mu_{\oo{t}}((W_2'(B)-A)P)=\mu_{\oo{t}}\otimes\mu_{\oo{t}}(D_BP),$$
with $D_A$ (resp. $D_B$) the non-commutative derivative with respect
to $A$ (resp. $B$) $\mu_A$ (resp. $\mu_B$) and  $W_i(z)=z^2/2+V_i(z)$.
Now, let $\mu_A$ (resp. $\mu_B$) be the spectral measure of the matrix $A$
(resp. $B$) then we shall obtain an algebraic equation for 
$H\mu_A(x)$ (resp. $H\mu_B(x)$) the Stieljes
transform of the limiting measure $\mu_A$ (resp. $\mu_B$)
given, for $x\in \CC\backslash \RR$ by:
$$H\mu_A(x)=\mu_{\oo{t}}(\frac{1}{ x-A})=\int \frac{1}{x-y}d\mu_A(y)$$

\begin{pr}\label{eyn}
Let for $x,y\in \CC\backslash \RR$,
$Y(x)=H\mu_A(x)-W_1'(x)$ and $X(y)=H\mu_B(x)-W_2'(x)$.
Then, there exists a polynomial function 
$$E(x,y)=\sum_{i,j=1}^{d-1} a_{ij}({\oo{t}}) x^i y^j$$
so that for all $x,y\in \CC\backslash \RR$
$$E(X(y),y)=0\qquad  E(x, Y(x))=0.$$
In particular, $\mu_A$ and $\mu_B$ 
are absolutely continuous with respect to Lebesgue 
measure, with Hilbert transform $H\mu_A$ and $H\mu_B$
so that $Y(x)=H\mu_A(x) -W_1'(x)$
satisfies the same algebraic equation with $x\in \RR$.
\end{pr}
\begin{dem}
Note that since we know
that $\mu_{\oo{t}}$ is compactly supported,
we can take in {\bf SD[$V_{Ising}$]}
Stieljes functions instead of polynomials $P$
since the latest are dense by Weirstrass theorem.

We take $P=P(A)=( x-A)^{-1}$ in the second equation 
in {\bf SD[$V_{Ising}$]}
to obtain:
$$\mu_{\oo{t}}\left(\frac{W_2'(B)}{ x-A}\right)=-1+xH\mu_A(x)$$
Then we use this in the first equation written with 
$$P(A,B)=\frac{1}{(x-A)}\frac{(W_2'(y)-W_2'(B))}{(y-B)}$$
to get
after some calculation
\begin{equation}\label{poilu}
U(x,y)(y-Y(x))=(Y(x)-W_1'(x))(x-W_2'(y))+1-Q(x,y)
\end{equation}
where
$$U(x,y)=\mu_{\oo{t}}\left(\frac{1}{(x-A)}\frac{W_2'(y)-W_2'(B)}{(y-B)}\right),$$
and
$$Q(x,y)=\mu_{\oo{t}}\left(\frac{W_1'(x)-W_1'(A)}{(x-A)}\frac{W_2'(y)
-W_2'(B)}{(y-B)}\right).$$
To obtain our algebraic equation, we simply define
$$E(x,y)=(Y(x)-W_1'(x))(x-W_2'(y))+1-Q(x,y)$$
and we obtain the famous ``Master-loop equation''
$$E(x,Y(x))=0$$
by taking $y=Y(x)$ in \eqref{poilu}.
In a symmetric way, we can show that if $X(y)=H\mu_B
(x)-W_2'(y)$ then we also have 
$E(X(y),y)=0$.
Note that $E$ is a polynomial function.
Hence, this shows that $Y(x)$, $X(y)$
and so the generating 
functions $H\mu_A(x)$ and $H\mu_B(y)$ are solution to an 
algebraic equation. However, this equation
still contains a certain
numbers of unknown; $\{\mu_{\oo{t}}(A^p B^q), p\le \mbox{deg}(V_1)-2, q
\le \mbox{deg}(V_2)-2\}.$ 
It is argued in physics 
that when ${\oo t}$ is small,
the supports of $\mu_A$ and $\mu_B$ should be connected
and therefore $(x,Y(x))$ and $(X(y),y)$
should then be genus zero curves. Then, these unknowns
should be determined by
the asymptotic behaviour
of $X(y)$ and $Y(x)$  at infinity
$$X(y)\simeq W_2'(y)-\frac{1}{y}(1+o(1)),\quad
Y(x)\simeq W_1'(x)-\frac{1}{x}(1+o(1)).
$$

Note in passing that, as solutions of an algebraic
equation, $H\mu_A$ and $H\mu_B(x)$ 
extends continuously (but in general not differentially)
 to the real line (eventually as an extended complex
number). As a consequence, $\mu_A$ and $\mu_B$
have densities with respect to the Lebesgue measure,
as the limits of the imaginary part
of the Stieljes transform 
on the real line.
\end{dem}

\subsection{Large deviations approach}
A large deviation approach was developed
in \cite{GCMP}, see also
 Matytsin \cite{matytsin}.
Again, we take $c=1$ up to rescaling
and denote $W_i(x)=x^2/2+V_i(x)$ for $i=1,2$.
 The main advantage of this
strategy is to be valid 
in the whole range of the
parameters. Otherwise, it should provide
the same type of information 
 than in the 
previous paragraph. Namely,

\begin{pr}
For any polynomials $V_1,V_2$ going to infinity
faster than $x^2$, $\mun_{A,B}$ 
converges almost surely 
towards $\mu_{\oo{t},1}=\mu_{\oo{t}}$ 
which is uniquely defined
by the Schwinger-Dyson's equations

\begin{equation}\label{masterloop}
\mu_{\oo{t}}\otimes\mu_{\oo{t}}(D_A P)=
\mu_{\oo{t}}((W_1'(A)- B )P),\qquad
\mu_{\oo{t}}\otimes\mu_{\oo{t}}(D_BP)=
\mu_{\oo{t}}((W_2'(B)- A )P)
\end{equation}
and by the fact that 
$\mu_{\oo{t}}|_A$ and $\mu_{\oo{t}}|_B$ (which
are the limits of $\mun_A$ and $\mun_B$ respectively)
are  the unique minimizers
of
$$S^{V_1,V_2}(\mu)=\mu_A(W_1
)+\mu_B( W_2)
-2^{-1}\int\int\log|x-y|d\mu^A(x)d\mu^A(y)$$
$$
-2^{-1}\int\int\log|x-y|d\mu^B(x)d\mu^B(y) 
+\frac{1}{2}\inf_{\rho,m}\{ 
\int_0^1 \int \frac{m_t(x)^2}{\rho_t(x)} dxdt
+\frac{\pi^2}{3}
\int_0^1 \int {\rho_t(x)^3}dxdt\}$$
where the inf is taken over $m,\rho$ so that $\mu_t(dx)=
\rho_t(x)dx\in \Pa(\RR)$, $\mu_0(x\in.)=
\mu_A
(x\in .)$, $
\mu_1(x\in.)=\mu_B
(x\in .)$,  and
$$\partial_t
\rho_t(x)+
\partial_x m_t(x)=0.$$
The infimum in $(\rho_., m_.)$
is taken along the solution to 
a complex Burgers equation; let $\Omega=\{ x\in \R, t\in (0,1):
\rho_t(x)>0\}$ and define on $\Omega$ $u_t(x)=\rho_t(x)^{-1} m_t(x)$
and 
 $f_t(x)=u_t(x)
+i\pi \rho_t(x)$. Then on $\Omega$,
$$\partial_t f_t(x)+ f_t(x)\partial_x f_t(x)=0.$$
Moreover, with $\mu_A=\mu_{\oo{t}}|_{A}$ and $\mu_B=\mu_{\oo{t}}|_{B}$,
for $\mu_A$-almost all $x$
\begin{equation}\label{ml1}
W_1'(x)- u_0(x) =H\mu_A(x),\quad\mu_A\,\mbox{a.s.,}\quad
 W_2'(x)+ u_1(x) =H\mu_B(x),\quad\mu_B\,\mbox{a.s.}
\end{equation}

\end{pr}
In comparison with the previous statements,
we note that the above results
hold for all $c$ and $V_1,V_2$,
and not only for small parameters.

\nn
\goodbreak
\begin{dem}
Most of the proof is contained 
 in \cite{GCMP} where the convergence 
of $\mun_A$, $\mun_B$ towards the unique minimizers
of $S^{V_1,V_2}$ was proved (see Theorem 3.3 in  \cite{GCMP}),
as well as the fact that the limit is compactly supported 
and that $\mu_{\oo{t}}$ satisfies (\ref{masterloop})
but for $P\in\Ca^m_{st}({\mathbb R})$
(see section 3.2.1, p. 555
and 558, in \cite{GCMP}). It clearly extends 
to polynomial functions since $\mu_{\oo{t}}$
is compactly supported as its marginals are.
The only point we stress here is that
this imply that  $\mu_{\oo{t}}$
is also uniquely determined.
Indeed, by proceeding by induction over
 the degree in $B$ of a monomial
 function $P$,
 we see that
 $$\tau( B P)=-\tau\otimes\tau(D_AP)+\tau(W_1'(A)P)$$
 defines uniquely all the moments $\tau(P(A,B))$
 from those of $\tau(Q(A))$.
 Note here that this is specific 
to the interaction under consideration;
in general the solutions of {\bf SD[V]} is not determined 
by their restriction to one variable.
\end{dem}

Using for instance 
the fact that if we let
 $g_t(x)=t f_t(x) +x$, the Wronskian
of $(f,g)$ is null, we find that
on each connected component
of $\Omega$, there exists an
analytic function $F$ so that 
$$t f_t(x) +x=F(f_t(x)).$$
In a small parameter region, it should easily 
be arguable that $\Omega$ is connected,
as it is when the parameters are null
(where the solution at time $t$  can be seen
to be a semi-circular variable
with variance $1-t+t^2$).
According to the previous section, 
we know that $f_t$ extends continuously to $t=0$ and $t=1$
since $\mu_A$ and $\mu_B$ have densities 
which yields 
\begin{equation}\label{riemann}
x=F(f_0(x))\quad f_1(y)+y=F(f_1(y))
\end{equation}
for all $x$ in the support of $\mu_A$ 
 and all $y$ in the support of $\mu_B$.
Noting that $f_0(x)=W_1'(x) -\overline{H\mu_A}(x)=-\overline{Y(x)}$,
$f_1(x)=-W_2'(x)+H\mu_B(x)=-X(x)$ 
it is tempting to hope that (\ref{riemann})
yields the same result that Property \ref{eyn},
namely that $(Y(x),x)$ and $(y,X(y))$ satisfy the
same algebraic equation. Our knowledge of this field
is much too limited to unable us to get this conclusion.

\nn{\bf Acknowledgments:}
We are extremely grateful
to B. Eynard and G. Schaeffer for many
comments which helped us to
compare the different mathematical
approaches to the enumeration of planar
maps. We also thank A. Okounkov
for many useful discussions.

\bibliography{map}
\bibliographystyle{acm}

\end{document}